\documentclass{tran-l}
\usepackage{amsmath,amssymb,amscd}
\usepackage[dvips]{graphics,epsfig}
\usepackage[all]{xy}

\newtheorem{theorem}{Theorem}[section]
\newtheorem{lemma}[theorem]{Lemma}
\newtheorem{prop}[theorem]{Proposition}
\newtheorem{cor}[theorem]{Corollary}

\theoremstyle{definition}
\newtheorem{definition}[theorem]{Definition}
\newtheorem{example}[theorem]{Example}

\theoremstyle{remark}
\newtheorem{remark}[theorem]{Remark}

\numberwithin{equation}{section}

\newcommand{\pa}{\partial}
\newcommand{\LL}{{\mathbb L}}
\newcommand{\LLL}{{\mathcal L}}
\newcommand{\Z}{{\mathbb Z}}
\newcommand{\HH}{{\mathbb H}}
\newcommand{\C}{{\mathbb C}\,}
\newcommand{\PP}{{\mathbb P}}
\newcommand{\Q}{{\mathbb Q}}

\newcommand{\I}{{\mathcal I}}
\newcommand{\ST}{{\mathcal S}}
\newcommand{\TT}{{\mathcal T}}
\newcommand{\transverse}{\pitchfork}

\begin{document}

\title{Homology manifold bordism}

\author{Heather Johnston}
\address{H.J.: Dept. of Mathematics\newline
\indent University of Massachusetts\newline
\indent Amherst, MA 01003\newline
\indent USA}
\email{johnston@math.umass.edu}
\thanks{This work was carried out in connection with the first named
author's EPSRC Visiting Fellowship in Edinburgh in August, 1997.}

\author{Andrew Ranicki}
\address{A.R.: Dept. of Mathematics and Statistics\newline
\indent University of Edinburgh\newline
\indent Edinburgh EH9 3JZ\newline
\indent Scotland, UK}
\email{aar@maths.ed.ac.uk}


\date{March 23, 1998 and, in revised form, September 3, 1999.}


\begin{abstract}
The Bryant-Ferry-Mio-Weinberger surgery exact sequence for compact
$ANR$ homology manifolds of dimension $\geq 6$ is used to obtain
transversality, splitting and bordism results for homology manifolds,
generalizing previous work of Johnston.

First, we establish homology manifold transversality for submanifolds
of dimension $\geq 7$: if $f:M \to P$ is a map from an $m$-dimensional
homology manifold $M$ to a space $P$, and $Q \subset P$ is a subspace
with a topological $q$-block bundle neighborhood, and $m-q \geq 7$, then
$f$ is homology manifold $s$-cobordant to a map which is transverse to $Q$,
with $f^{-1}(Q) \subset M$ an $(m-q)$-dimensional homology submanifold.

Second, we obtain a codimension $q$ splitting obstruction $s_Q(f) 
\in LS_{m-q}(\Phi)$ in the Wall $LS$-group for a simple homotopy 
equivalence $f:M \to P$ from an $m$-dimensional homology manifold 
$M$ to an $m$-dimensional Poincar\'e space $P$ with a codimension 
$q$ Poincar\'e subspace $Q \subset P$ with a topological normal bundle,  
such that $s_Q(f)=0$ if (and for $m-q \geq 7$ only if) $f$ splits at $Q$ 
up to homology manifold $s$-cobordism.

Third, we obtain the multiplicative structure of the homology manifold
bordism groups $\Omega^H_*\cong\Omega^{TOP}_*[L_0(\Z)]$.

\end{abstract}

\maketitle

\section*{Contents}
\begin{itemize}
\item[~] Introduction
\item[1.] Homology manifold surgery and bordism
\item[2.] Homology submanifold transversality up to bordism
\item[3.] Codimension $q$ splitting of homology manifolds
\item[4.] Homology submanifold transversality up to $s$-cobordism
\item[5.] Dual transversality for homology manifolds
\item[6.] The dual transversality obstruction
\item[7.] The multiplicative structure of homology manifold bordism
\item[~] References
\end{itemize}

\section*{Introduction}

Homology manifolds are spaces with the local homology properties of
topological manifolds, but not necessarily their geometric properties
such as transversality.  The results of Johnston \cite{Jo} on the bordism
and transversality properties of high-dimensional homology manifolds
are extended here using the methods of surgery theory.  The extent to which
transversality holds in a homology manifold is a measure of how close it is
to being a topological manifold. It is not possible to investigate
transversality in homology manifolds by direct geometric methods --
as in \cite{Jo} we employ bordism and surgery instead.

\medskip

We start with a brief recollection of transversality for differentiable
manifolds.  Suppose that $P^n$ is an $n$-dimensional differentiable
manifold and $Q^{n-q} \subset P^n$ is a codimension $q$ submanifold
with a $q$-plane normal bundle
$$\nu_{Q \subset P}~:~Q \to BO(q)~.$$
A smooth map $f:M \to P$ from an $m$-dimensional differentiable
manifold $M^m$ is transverse to $Q \subset P$ if the inverse image of
$Q$ is a codimension $q$ submanifold
$$N^{m-q}~=~f^{-1}(Q) \subset M^m$$
with normal $q$-plane bundle the pullback of $\nu_{Q\subset P}$
along $g=f\vert:N \to Q$
$$\nu_{N \subset M}~=~g^*\nu_{Q \subset P}:N \to BO(q)~.$$
The classical result of Thom is that every map $f:M^m \to P^n$ is
homotopic (by a small homotopy) to a smooth map which is transverse to
$Q \subset P$.  This result was proved by direct analytic methods.

\medskip

Topological manifolds also have transversality, but the proof is very
indirect, relying heavily on surgery theory -- see Kirby and Siebenmann
\cite{KS} (III,\S 1), Marin \cite{M} and Quinn \cite{Q3}.  Instead of
vector bundles it is necessary to work with normal microbundles,
although we shall use the formulation in terms of the topological block
bundles of Rourke and Sanderson \cite{RS}.

\medskip

The essential aspect of transversality is that a submanifold has a
nice normal (vector or block) bundle, as formalized in the
following definition.

\medskip

\noindent{\bf Definition}
A {\it codimension $q$ bundle subspace} $(Q,R,\xi)$ ($q\geq 1$) of a
space $P$ is a subspace $Q \subset P$ together with a topological
$q$-block bundle
$$(D^q,S^{q-1}) \to (E(\xi),S(\xi)) \to Q$$
such that
$$P~=~E(\xi)\cup_{S(\xi)}R$$
where $R=\overline{P\backslash E(\xi)}$. When $R$ and $\xi$ are clear
we say that $Q$ is a codimension $q$ bundle subspace of $P$.
\medskip

Topological $q$-block bundles over a space $Q$ are classified by the
homotopy classes of maps from $Q$ to a classifying space
$B\widetilde{TOP}(q)$, so we write such a bundle $\xi$ as a map
$$\xi~:~Q \to B\widetilde{TOP}(q)~.$$
If $P$ is an $n$-dimensional topological manifold and $Q \subset P$
is a triangulable locally flat codimension $q$ submanifold with 
$n-q \geq 5$ or $q \leq 2$ then $Q$ is a codimension $q$ bundle 
subspace of $P$ with
$$\xi~=~\nu_{Q\subset P}~:~Q \to B\widetilde{TOP}(q)$$
a normal topological $q$-block bundle, by Theorem 4.9 of Rourke and
Sanderson \cite{RS}.  (Hughes, Taylor and Williams \cite{HTW} obtained
a topological regular neighborhood theorem for arbitrary locally flat
submanifolds in a manifold of dimension $\geq 5$, in which the
neighborhood is the mapping cylinder of a manifold approximate
fibration).

\medskip

In the applications of codimension $q$ bundle subspaces $(Q,R,\xi)\subset P$
we shall only be concerned with the case when $P$ is a finite $CW$
complex and $Q,R \subset P$ are subcomplexes.

\medskip

\noindent{\bf Definition}  (Submanifold transversality)\\
Let $P$ be a space with a codimension $q$ bundle subspace $(Q,R,\xi)$.\\
{\rm (i)} A map $f:M \to P$ from an $m$-dimensional manifold $M$
is {\it transverse} to $Q\subset P$ if the inverse image of $Q$
$$N^{m-q}~=~f^{-1}(Q) \subset M^m$$
is a locally flat codimension $q$ submanifold with the pullback normal bundle.\\
{\rm (ii)} A map $f:M \to P$ is {\it $s$-transverse} to $Q \subset P$
if it is $s$-cobordant to a transverse map.
\medskip

Of course, the submanifolds of the manifold $M$ and bundles in the above
definitions are understood to be in the same category as $M$ itself.
For simplicity, we shall only be considering compact oriented manifolds.
\medskip

An {\it $m$-dimensional homology manifold} $M$ is a finite-dimensional
$ANR$ such that for each $x \in M$
$$H_r(M,M\backslash \{x\})~=~
\begin{cases}
\Z&\text{if $r=m$} \\
0&\text{if $r \neq m$}~.
\end{cases}
$$
\smallskip

We shall say that an $m$-dimensional homology manifold $M$ has {\it
codimension $q$ $s$-transversality} if every map $f:M \to P$ is
$s$-transverse to every codimension $q$ bundle subspace $Q \subset P$.  
(It is unknown if the analogue of the topological $s$-cobordism theorem 
holds for homology manifolds).

\medskip

An $m$-dimensional homology manifold $M$ is {\it resolvable} if there
exists a $CE$ map $h:M_{TOP} \to M$ from an $m$-dimensional topological
manifold $M_{TOP}$. (Roughly speaking, a $CE$ map is a map with contractible
point inverses). Resolvable homology manifolds have codimension $q$
$s$-transversality for all $q \geq 1$\ :
if $f:M \to P$ is a map from a resolvable $m$-dimensional
homology manifold and $Q \subset P$ is a codimension $q$ bundle subspace,
then the mapping cylinder of $h$ is a homology manifold $s$-cobordism
$$(g;f,f_{TOP})~:~(W;M,M_{TOP}) \to P$$
from $f$ to a map $f_{TOP}:M_{TOP} \to P$ which can be made (topologically)
transverse to $Q \subset P$.

\medskip

Quinn \cite{Q1} used controlled surgery to prove that for $m\geq 6$
an $m$-dimensional homology manifold $M$ with
codimension $m$ $s$-transversality is resolvable.
The {\it resolution obstruction} of Quinn \cite{Q2}
$$i(M) \in H_m(M;L_0(\Z))$$
is such that $i(M)=0$ if (and for $m \geq 6$ only if) $M$ is
resolvable; for connected $M$ the obstruction takes values in
$H_m(M;L_0(\Z))=\Z$.  The invariant $i(M)$ is the obstruction to a degree 1
map $f:M^m \to S^m$ being $s$-transverse to some point $* \in S^m$. 
Bryant, Ferry, Mio and Weinberger \cite{BFMW} constructed exotic
homology manifolds $M^m$ in dimensions $m\geq 6$ which are not
resolvable, and initiated the surgery classification theory for
high-dimensional homology manifolds up to $s$-cobordism.  In Chapter 1
we shall modify the construction of \cite{BFMW} to obtain a connected
homology manifold $M=N_I$ with prescribed resolution obstruction $I \in
L_0(\Z)$, starting with any connected $m$-dimensional topological
manifold $N$ ($m \geq 6$).  This homology manifold is not homotopy
equivalent to $N$, but it is in a prescribed homology manifold normal
bordism class of $N$.

\medskip

The first named author used the theory of \cite{BFMW} to prove that
$m$-dimensional homology manifolds have codimension $q$ $s$-transversality
and splitting in the following cases.

\medskip

\noindent{\bf Theorem}
{\rm (Homology manifold $\pi$-$\pi$ $s$-transversality, Johnston \cite{Jo})}\\
{\it Let $f:M \to P$ be a map from an $m$-dimensional
homology manifold $M$ to a space $P$ with a codimension $q$ bundle subspace $(Q,R,\xi)$,
with $m-q \geq 6$.\\
{\rm (i)} If $q=1$, $\xi$ is trivial, and $R=R_1 \sqcup R_2$ is
disconnected with $\pi_1(Q)\cong \pi_1(R_1)$, then $f$ is
$s$-transverse to $Q \subset P$.\\
{\rm (ii)} If $q\geq 3$ then $f$ is $s$-transverse to $Q \subset P$.}

\medskip

\noindent{\bf Definition} (i) A {\it codimension $q$ Poincar\'e bundle
subspace} $(Q,R,\xi)$ of an $m$-dimen\-sional Poincar\'e space $P$ is a
codimension $q$ bundle subspace such that $Q$ is an $(m-q)$-dimensional
Poincar\'e space and $(R,S(\xi))$ is an $m$-dimensional Poincar\'e
pair, where $S(\xi)$ is the total space of the $S^{q-1}$-bundle of
$\xi$ over $P$.\\
(ii) Let $P,Q,R,\xi$ be as in (i).  A simple homotopy equivalence $f:M
\to P$ from an $m$-dimensional homology manifold $M$ {\it splits at} $Q
\subset P$ if $f$ is $s$-cobordant to a simple homotopy equivalence
(also denoted by $f$) which is transverse to $Q \subset P$ and such
that the restrictions $f\vert : f^{-1}(Q) \to Q$, $f\vert : f^{-1}(R)
\to R$ are also simple homotopy equivalences.  
\medskip

\noindent{\bf Theorem} 
{\rm (Homology manifold splitting, Johnston \cite{Jo})}\\
{\it Let $f:M \to P$ be a simple homotopy equivalence from an $m$-dimensional
homology manifold $M$ to an $m$-dimensional Poincar\'e
space $P$ with a codimension $q$ Poincar\'e bundle subspace $(Q,R,\xi)$,
with $m-q \geq 6$. \\
{\rm (i)} {\rm (Codimension 1 $\pi$-$\pi$ splitting)} 
If~$q=1$, $\xi$ is trivial, and $R=R_1 \sqcup R_2$ is disconnected with 
$\pi_1(Q)\cong \pi_1(R_1)$, then $f$ splits at $Q \subset P$.\\
{\rm (ii)} {\rm (Browder splitting)} If~$q\geq 3$ then $f$ splits at $Q \subset P$
if and only if the restriction $g=f\vert:f^{-1}(Q) \to Q$ has
surgery obstruction $\sigma_*(g)=0\in L_{m-q}(\Z[\pi_1(Q)])$.}

\medskip

In Chapters 2,3,4 we shall use the theory of \cite{BFMW} to obtain
even stronger results on homology manifold transversality and
codimension $q$ splitting.  The results of this paper require a
slightly higher dimension hypothesis $m-q\geq 7$.  The extra dimension
is needed to apply codimension 1 $\pi$-$\pi$ splitting, (i) above, to
a homology manifold of dimension $m-q$. Thus we require $m-q-1\geq 6$
or $m-q\geq 7$.

\medskip

Wall \cite{Wa} (Chapter 11) obtained a codimension $q$ splitting
obstruction $s_Q(f) \in LS_{m-q}(\Phi)$ for a simple homotopy
equivalence $f:M \to P$ from an $m$-dimensional topological manifold
$M$ to an $m$-dimensional Poincar\'e space $P$ with a codimension
$q\geq 1$ Poincar\'e bundle subspace $(Q,R,\xi)$, such that $s_Q(f)=0$
if (and for $m-q \geq 5$ only if) $f$ splits at $Q\subset P$.  Our
first main result obtains the analogous obstruction for the
codimension $q$ splitting of homology manifolds.  (The full statement
will be given in Theorem \ref{qsplit}.)
\medskip

\begin{theorem} \label{fifteen}
A simple homotopy equivalence $f:M \to P=E(\xi)\cup_{S(\xi)}R$
from an $m$-dimensional homology manifold $M$ to an $m$-dimensional Poincar\'e
space $P$ with a codimension $q$ Poincar\'e bundle subspace $(Q,R,\xi)$
has  a codimension $q$ splitting obstruction
$$s_Q(f) \in LS_{m-q}(\Phi)$$
such that $s_Q(f)=0$ if {\rm (}and for $m-q \geq 7$ only if\/{\rm )}
$f$ splits at $Q \subset P$.
\end{theorem}

Our second main result establishes homology manifold
$s$-transversality in the case $m-q\geq 7$, generalizing the homology
manifold $\pi$-$\pi$ $s$-transversality theorem of \cite{Jo}. This
result appears as Theorem \ref{homtra} below.

\begin{theorem} \label{one}
Let $f:M \to P=E(\xi)\cup_{S(\xi)}R$ be a map from an $m$-dimensional
homology manifold $M$ to a space $P$ with a codimension $q$ bundle subspace $(Q,R,\xi)$.
If $m-q\geq 7$ then $f$ is $s$-transverse to $Q \subset P$.
\end{theorem}

In Chapters 5,6
we consider $s$-transversality for a map $f:M \to P$ from a homology manifold
$M$ to the polyhedron $P=\vert K\vert$ of a (finite) simplicial complex $K$.
Instead of seeking $s$-transversality to just one codimension
$q$ bundle subspace $Q \subset P$ we consider $s$-transversality to all the
dual cells $\vert D(\sigma,K)\vert \subset P$ ($\sigma \in K$) at
once, following the work of Cohen \cite{Co} on $PL$ manifold transversality.

\medskip

The {\it dual cells} of a simplicial complex $K$ are the subcomplexes of $K'$
$$D(\sigma,K)~=~
\{\widehat \sigma_0 \widehat \sigma_1 \dots \widehat \sigma_p \in
K' \, \vert \,  \sigma  \le  \sigma_0 < \sigma_1 < \dots < \sigma_p \}~~
(\sigma \in K)~.$$

The boundary of a dual cell is the subcomplex
$$\begin{array}{rl}
\partial D(\sigma ,K)\!\!&=~ \bigcup_{\tau > \sigma} D(\tau ,K)\vspace*{2mm}\\
&=~\{\widehat \sigma_0 \widehat \sigma_1 \dots \widehat \sigma_p \in
K' \, \vert \,  \sigma  <  \sigma_0 < \sigma_1 < \dots < \sigma_p \}
\subset D(\sigma,K)~.
\end{array}$$
\smallskip

\noindent{\bf Definition} (Dual transversality)\\
(i) A map $f: M \to \vert K\vert$ from an $m$-dimensional manifold $M$
(in some manifold category)
is {\it dual transverse} if the inverse images of the dual cells are
codimension $\vert \sigma \vert$ submanifolds
$$M(\sigma)^{m-\vert\sigma \vert}~=~f^{-1}(D(\sigma,K)) \subseteq M^m$$
with boundary
$$\partial M(\sigma)~=~f^{-1}(\partial D(\sigma,K)) ~=~
\bigcup\limits_{\tau >\sigma}M(\tau)~.$$
(ii) An $m$-dimensional manifold $M$ has {\it dual $s$-transversality}
if every map $f:M \to \vert K \vert$ is $s$-cobordant to a dual
transverse map.

\medskip

Dual transversality implies submanifold transversality\ : if
$$f~:~M^m \to P~=~\vert K \vert ~=~ E(\xi) \cup_{S(\xi)}R$$
is dual transverse then $f$ is transverse to every polyhedral
codimension $q$ bundle subspace $Q\subset P$.  $PL$ manifolds and $PL$
homology manifolds $M$ have dual transversality, with every simplicial
map $f:M \to K'$ dual transverse -- in this case each inverse image
$f^{-1}(D(\sigma,K))\subset M$ ($\sigma \in K$) is automatically a $PL$
submanifold of codimension $\vert \sigma \vert$ (Cohen \cite{Co}), so
there is no need to use $s$-cobordisms.  Topological manifolds $M$ have
dual transversality by the work of Kirby-Siebenmann \cite{KS} and Quinn
\cite{Q3}, with every map $f:M \to \vert K\vert$ homotopic to a dual
transverse map.

\medskip

The $s$-transversality result of Theorem \ref{one} can be applied
inductively to obtain dual $s$-transversality for a map $f:M^m \to
\vert K \vert$ in the case when the inverse images
$f^{-1}(D(\sigma,K))$ ($\sigma \in K$) are required to be homology
manifolds of dimensions $m - \vert\sigma \vert \geq 7$.

\begin{cor} If $f:M \to \vert K \vert$ is a map from an $m$-dimensional
homology manifold $M$ to the polyhedron of a $k$-dimensional simplicial complex
$K$ with $m-k \geq 7$ then $f$ is dual $s$-transverse.
\end{cor}

On the other hand, if $m-k <7$, dual transversality may be obstructed. 
Consider for example the case of $f:M\to \vert K \vert$ for $m=k$, with
$M$ connected.  In this case the resolution obstruction is easily shown
to be an obstruction to dual transversality.  If $f:M\to \vert K \vert$
is dual transverse, then for some $m$-simplex $\sigma^m\in K$ and a
sequence of faces $\sigma^0 < \sigma^1 < \dots < \sigma^m$
the inverse images $M(\sigma^j)^{m-j}=f^{-1}(D(\sigma^j,K))$ are non-empty 
codimension $j$ homology submanifolds of $M$ with
$$M(\sigma^m)^0 \subset \partial M(\sigma^{m-1}) \subset M(\sigma^{m-1})^1
\subset \dots \subset M(\sigma^0)^m \subset M~.$$
Quinn (\cite{Q2}, 1.1) proved that the resolution obstruction $i(X)$ is 
such that\,: 
\begin{itemize}
\item[(i)] $i(U)=i(X)$ for any open subset $U \subset X$ of a homology 
manifold $X$,
\item[(ii)] $i(X)=i(\partial X)$ for any 
connected homology manifold $X$ with non-empty boundary $\partial X$.
\end{itemize}
It follows that
$$i(M(\sigma^m)^0)~=~i(M(\sigma^{m-1})^1)~=~\dots~=~i(M(\sigma^0)^m)~=~i(M)~.$$  
Now $M(\sigma^m)$ is a 0-dimensional homology manifold, which is a
(finite) union of points, so that $i(M)=i(M(\sigma^m))=0$ and $M$ is
resolvable.  In Chapter 6 we use the algebraic topology of homology
manifold bordism to prove a strong generalization of this result. 
Although the resolution obstruction continues to play a key role, the
general form of the dual transversality obstruction is more complicated
algebraically.

\medskip

For any space $X$ let $\Omega^H_m(X)$ be the bordism group of maps
$M \to X$ from $m$-dimensional homology manifolds.
An $m$-dimensional topological manifold is an $m$-dimensional homology
manifold, so there are evident forgetful maps
$$\Omega^{TOP}_m(X) \to \Omega^H_m(X)~.$$

For any simplicial complex $K$, let $\Omega^{H,\transverse}_m(K)$
be the bordism group of dual transverse maps
$M \to \vert K \vert$ from $m$-dimensional homology manifolds.
Forgetting dual transversality gives maps
$$A^H~:~\Omega^{H,\transverse}_*(K) \to \Omega^H_*(K)~.$$
The extent to which dual transversality holds for homology manifolds up
to bordism is measured by the extent to which the maps $A^H$ are
isomorphisms.  
\medskip

Our third main result relates the obstruction to homology manifold dual
$s$-transversality to the resolution obstruction, and
identifies the fibre of $A^H$ with a generalized homology theory.
The full statement will be given in Theorem \ref{two},
including the following result\,:

\begin{theorem} For $m \geq 6$ the $K$-transverse homology manifold
bordism groups $\Omega^{H,\transverse}_m(K)$ are related to the
homology manifold bordism groups $\Omega^H_m(K)$ by an exact sequence
$$\dots \to\Omega^{H,\transverse}_m(K) \to \Omega^H_m(K)
\to H_m(K;\LLL_{\bullet})\to \Omega^{H,\transverse}_{m-1}(K) \to \dots$$
with $\LLL_{\bullet}$ a spectrum such that
$\pi_0(\LLL_{\bullet})=\Z[L_0(\Z)\backslash \{0\}]$ and
$\pi_m(\LLL_{\bullet})=0$ for $m\geq 6$.
\end{theorem}

Ferry and Pedersen \cite{FP} showed that the Spivak normal fibration $\nu_M$
of an $m$-dimensional homology manifold $M$ admits a canonical $TOP$
reduction $\widetilde{\nu}_M$, so that there is a canonical bordism class
of normal maps $M_{TOP} \to M$ from topological manifolds. The surgery
obstruction $\sigma_*(M_{TOP} \to M) \in L_m(\Z[\pi_1(M)])$ is
determined by the resolution obstruction $i(M) \in H_m(M;L_0(\Z))$.

\medskip

For any abelian group $A$ let $A[L_0(\Z)]$ be the abelian group
of finite linear combinations
$$\sum\limits_{i \in L_0(\Z)}a_i[i] ~~(a_i \in A)~,$$
that is the direct sum of $L_0(\Z)$ copies of $A$
$$A[L_0(\Z)]~=~\Z[L_0(\Z)]\otimes_{\Z}A~=~\bigoplus_{L_0(\Z)}A~.$$
\medskip

\noindent{\bf Theorem} 
{\rm (The additive structure of homology manifold bordism, Johnston \cite{Jo})}\\
{\it For any simplicial complex $K$ the map of bordism groups
$$\Omega^H_m(K) \to \Omega_m^{TOP}(K)[L_0(\Z)]~;~
(M,f) \mapsto (M_{TOP},f_{TOP})[i(M)]$$
is an isomorphism for $m \geq 6$, with}
$$f_{TOP}~:~M_{TOP} \to M \stackrel{f}{\to} \vert K\vert~.$$

\medskip
In Chapter 7 we shall analyze the multiplicative structure
on $\Omega_m^{TOP}(K)[L_0(\Z)]$ which corresponds to the
cartesian product of homology manifolds under the
isomorphism $\Omega^H_m(K) \cong \Omega_m^{TOP}(K)[L_0(\Z)]$.
\medskip

We should like to thank the referee and Bruce Hughes for 
helpful comments on the manuscript.

\section{Homology manifold surgery and bordism}

We review and extend the surgery theory of $ANR$ homology manifolds.
\medskip

An (oriented) {\it simple $m$-dimensional Poincar\'e space} $X$ is a compact
$ANR$ with a fundamental class $[X] \in H_m(X)$ such that the
chain map
$$[X] \cap - ~:~ C(\widetilde{X})^{m-*}~=~{\rm Hom}_{\Z[\pi_1(X)]}(C(\widetilde{X}),
\Z[\pi_1(X)])^{m-*} \to C(\widetilde{X})$$
is a simple $\Z[\pi_1(X)]$-module chain equivalence inducing isomorphisms
$$[X]\cap -~:~H^{m-*}(\widetilde{X})~\cong~H_*(\widetilde{X})~,$$
with $\widetilde{X}$ the universal cover of $X$ and
$H^*(\widetilde{X})\,\equiv\,H_{-*}(C(\widetilde{X})^{-*})$.
(A compact $ANR$ has a preferred simple homotopy type by the work of
Chapman \cite{Ch}). In particular, an $m$-dimensional homology manifold is
an $m$-dimensional Poincar\'e space.

\medskip

The {\it manifold structure set} $\ST(X)$ of a simple $m$-dimensional
Poincar\'e space $X$ is the set of equivalence classes of pairs $(M,h)$
with $M$ an $m$-dimensional topological manifold and $h:M \to X$ a
simple homotopy equivalence, with
$(M_1,h_1) \allowbreak \simeq
(M_2,h_2)$ if there exists an $s$-cobordism $(W;M_1,M_2)$ with a simple
homotopy equivalence of the type
$$(g;h_1,h_2)~:~(W;M_1,M_2) \to X \times ([0,1];\{0\},\{1\})~.$$
The {\it normal invariant set} $\TT(X)$ is the bordism set of degree 1
normal maps $(f,b):M \to X$ from topological manifolds, with $b:\nu_M
\to \widetilde{\nu}_X$ a bundle map from the normal bundle of $M$ to a
$TOP$ reduction of the Spivak normal fibration of $X$.  For $m\geq 5$
the Browder-Novikov-Sullivan-Wall surgery theory for topological
manifolds gives the surgery exact sequence
$$\dots \to L_{m+1}(\Z[\pi_1(X)]) \stackrel{\partial}{\to}
\ST(X) \stackrel{\eta}{\to} \TT(X)
\stackrel{\theta}{\to} L_m(\Z[\pi_1(X)])$$
(Wall \cite{Wa}, Chapter 10).  In general, it is possible for $\TT(X)$ and
$\ST(X)$ to be empty: the theory involves a primary topological
$K$-theory obstruction for deciding if $\TT(X)$ is non-empty and a
secondary algebraic $L$-theory obstruction for deciding if $\ST(X)$ is
non-empty.  More precisely, $\TT(X)$ is non-empty if and only if the
Spivak normal fibration $\nu_X:X \to BG(k)$ ($k$ large) admits a $TOP$
reduction $\widetilde{\nu}_X:X \to B\widetilde{TOP}(k)$, corresponding
by the Browder-Novikov transversality construction on a degree 1 map
$\rho:S^{m+k} \to T(\widetilde\nu_X)$ to a normal map
$$(f,b)~=~\rho\vert~ :~ M~ =~\rho^{-1}(X) \to X$$
from a topological manifold. A choice of $\widetilde{\nu}_X$ determines
a bijection
$$\TT(X)~\cong~ [X,G/TOP]~.$$
The algebraic surgery exact sequence of Ranicki \cite{R3}
$$\dots \to L_{m+1}(\Z[\pi_1(X)]) \to \ST_{m+1}(X)\to
H_m(X;\LL_\bullet) \stackrel{A}{\to} L_m(\Z[\pi_1(X)])\to \dots$$
is defined for any space $X$, with $A$ the assembly map on the
generalized homology group of the 1-connective quadratic $L$-theory
spectrum $\LL_\bullet=\LL_{\bullet}\langle 1 \rangle (\Z)$ of $\Z$,
with 0th space
$$\LL_0~ \simeq~ G/TOP$$
and homotopy groups the simply-connected surgery obstruction groups
$$\pi_m(\LL_{\bullet})~=~L_m(\Z)~=~
\begin{cases}
\Z&\text{if $m \equiv 0$ (mod 4)}\\
0&\text{if $m \equiv 1$ (mod 4)}\\
\Z_2&\text{if $m \equiv 2$ (mod 4)}\\
0&\text{if $m \equiv 3$ (mod 4)~.}
\end{cases}$$

The total surgery obstruction $s(X) \in \ST_m(X)$ of a simple
$m$-dimensional geometric Poincar\'e complex is such that $s(X)=0$ if
(and for $m \geq 5$ only if) $X$ is simple homotopy equivalent to an
$m$-dimensional topological manifold.  The surgery exact sequence of an
$m$-dimensional topological manifold $M$ is isomorphic to the
corresponding portion of the algebraic surgery exact sequence, with
$$\ST(M)~=~\ST_{m+1}(M)~~,~~\TT(M)~=~[M,G/TOP]~=~H_m(M;\LL_{\bullet})~.$$

\indent 
The surgery theory of topological manifolds was extended to
homology manifolds in Quinn \cite{Q1},\cite{Q2} and 
Bryant, Ferry, Mio and Weinberger \cite{BFMW}, using
the 4-periodic algebraic surgery exact sequence of Ranicki \cite{R3}
(Chapter 25)
$$\dots \to L_{m+1}(\Z[\pi_1(X)]) \to \overline{\ST}_{m+1}(X)\to
H_m(X;\overline{\LL}_\bullet) \stackrel{\overline{A}}{\to} 
L_m(\Z[\pi_1(X)])\to \dots~.$$
This sequence is defined for any space $X$, with $\overline{A}$ the assembly 
map on the generalized homology group of the 0-connective quadratic $L$-theory
spectrum $\overline{\LL}_\bullet=\overline{\LL}_{\bullet}\langle 0
\rangle (\Z)$ of $\Z$, with 0th space
$$\overline{\LL}_0~ \simeq~ G/TOP\times L_0(\Z)~.$$
The 4-periodic total surgery obstruction $\overline{s}(X) \in
\overline{\ST}_m(X)$ of a simple $m$-dimensional geometric Poincar\'e
complex $X$ is such that $\overline{s}(X)=0$ if (and for $m \geq 6$
only if) $X$ is simple homotopy equivalent to an $m$-dimensional
homology manifold, by \cite{BFMW}.  The {\it homology manifold
structure set} $\ST^H(X)$ of a simple $m$-dimensional Poincar\'e space
$X$ is defined in the same way as $\ST(X)$ but using homology
manifolds.  The surgery exact sequence of an $m$-dimensional homology
manifold $M$ is isomorphic to the corresponding portion of the
4-periodic algebraic surgery exact sequence, with
$$\ST^H(M)~=~\overline{\ST}_{m+1}(M)~.$$ 
\indent The essential difference between surgery on homology manifolds
and topological manifolds is that there is no Browder-Novikov
transversality allowing the construction of normal maps from homology
manifolds.  Thus, the surgery exact sequence of \cite{BFMW} does not
follow Wall \cite{Wa} in relating homology manifold structures and
normal invariants.  Rather, the homology manifold surgery exact
sequence of \cite{BFMW} follows the stratified surgery exact sequence
of Weinberger \cite{Wei} in that it relates two purely algebraically
defined groups with the geometrically defined structure set.  Despite
the fact that \cite{BFMW} does not define homology manifold normal
invariants, one can define homology manifold normal invariants,
$\TT(X)$ similar to normal invariants of topological manifolds. 
\medskip

The {\it homology normal invariant set} $\TT^H(X)$ is the bordism set
of degree 1 normal maps $(f,b):M \to X$ from connected homology
manifolds, with $b:\widetilde{\nu}_M \to \widetilde{\nu}_X$ a map from
the canonical $TOP$ reduction (\cite{FP}) of the Spivak normal
fibration of $M$ to a $TOP$ reduction of the Spivak normal fibration of
$X$.  It is still the case that $\TT^H(X)$ is non-empty if and only if
$\nu_X$ is $TOP$ reducible, but now it is necessary to also keep track
of the resolution obstruction, and the homology manifolds have to be
constructed using controlled topology.  The following theorem allows us
to use the geometric interpretation $\TT^H(X)$ in the homology manifold
surgery exact sequence of \cite{BFMW}.

\begin{theorem} \label{tth}{\rm (Johnston \cite{Jo})}\\
Let $m\geq 6$. For a connected simple $m$-dimensional Poincar\'e space
$X$ the function
$$\begin{array}{l}
\TT^H(X) \to \TT(X) \times L_0(\Z)~;\\
((f,b):M \to X) \mapsto ((f_{TOP},b_{TOP}):
M_{TOP} \to M \stackrel{(f,b)}{\to} X,i(M))
\end{array}
$$
is a natural bijection, with $M_{TOP} \to M$ the topological degree 1
normal map determined by the canonical $TOP$ reduction.  A choice of
$\widetilde{\nu}_X$ determines a bijection
$$\TT^H(X)~\cong~ [X,G/TOP\times L_0(\Z)]~.$$
\end{theorem}
Actually \cite{Jo} (5.2) is for $m\geq 7$, but we can improve to $m\geq 6$
by a slight variation of the proof as described below.
\medskip

Given the above theorem, the homology manifold surgery exact sequence
of \cite{BFMW} 
$$\dots \to L_{m+1}(\Z[\pi_1(X)]) \stackrel{\partial}{\to}
\ST^H(X) \stackrel{\eta}{\to} [X,G/TOP\times L_0(\Z)]
\stackrel{\theta}{\to} L_m(\Z[\pi_1(X)])$$
becomes: 
\begin{theorem}\label{eleven} 
{\rm (Bryant, Mio, Ferry and Weinberger \cite{BFMW})}\\
The homology manifold structure set $\ST^H(X)$ fits into the exact sequence 
$$\dots \to L_{m+1}(\Z[\pi_1(X)]) \stackrel{\partial}{\to}
\ST^H(X) \stackrel{\eta}{\to} \TT^H(X)
\stackrel{\theta}{\to} L_m(\Z[\pi_1(X)])~.$$
In particular, $\ST^H(X)$ is non-empty if and only if there exists a
degree 1 normal map $(f,b):M \to X$ from a homology manifold $M$ with
surgery obstruction $0 \in L_m(\Z[\pi_1(X)])$.
\end{theorem}

If $M$ is an $m$-dimensional homology manifold the canonical bordism
class of topological normal maps $M_{TOP} \to M$ determined by the
canonical $TOP$ reduction $\widetilde{\nu}_M$ (\cite{FP}) of the Spivak
normal fibration has surgery obstruction
$$\sigma_*(M_{TOP} \to M)~=~\overline{A}(-i(M)) \in L_m(\Z[\pi_1(M)])$$
the image (up to sign) of the resolution obstruction $i(M) \in H_m(M;L_0(\Z))$
under the composite
$$H_m(M;L_0(\Z)) \subseteq H_m(M;\overline{\LL}_{\bullet})
\stackrel{\overline{A}}{\to} L_m(\Z[\pi_1(M)])~.$$

\begin{remark} The algebraic surgery exact sequence has a geometric
interpretation: the map $\eta:\ST^H(X)\to \TT^H(X)$ is
given by the forgetful map and the map
$\partial:L_{m+1}(\Z[\pi_1(X)]) \to \ST^H(X)$ is given by 
``Wall Realization'', i.e. the homology manifold analogue of 
the constructions of Chapters 5 and 6 of \cite{Wa} realizing 
the elements of the $L$-groups as the surgery obstructions of normal
maps obtained by non-simply connected plumbing.
\end{remark}

To see that $\eta$ is the forgetful map take $f:Y\to X$ and compare
$\eta(f)$ with $(f:Y\to X)\in \TT^H(X)$.  By \ref{tth} above we may consider
instead their images in $\TT(X)\times L_0(\Z)$.  By definition of the
canonical $TOP$-reduction of the Spivak normal fibration $\nu_Y$
(Ferry and Pedersen \cite{FP}), $\eta(f)$ is given by $f\circ
f_{TOP}:Y_{TOP}\to Y\to X, I(Y)$, i.e.  it agrees with the image of
$f:Y\to X$.  
\medskip

``Wall Realization'' does not (yet) have an obvious geometric
interpretation for homology manifolds.  There are no homology manifold
``handles'' with nice attaching maps.  Nevertheless, one can obtain the
following theorem.

\begin{theorem} {\rm (Johnston \cite{Jo})}\\
Let $m \geq 6$.
For an $m$-dimensional Poincar\'e space $P$ with a specified homotopy
equivalence $h:X\to P$, with $X$ a homology manifold, and for any
element $\sigma \in L_{m+1}(\Z[\pi_1(P)])$, the image $\partial(\sigma)\in
\ST^H(P)$ under the map $\partial:L_{m+1}(\Z[\pi_1(P)])\to \ST^H(P)$ in the
surgery exact sequence of \cite{BFMW} has a representative $g:Y\to P$
such that there exists a homology manifold bordism $r:W\to
P\times[0,1]$ with $r|\pa W=g\amalg h$.
\end{theorem}

This theorem follows from \ref{tth} and the surgery exact sequence of
\cite{BFMW} for $\ST^H(P\times[0,1], P\times\{0,1\})$ relative to the
given structures $(g,h)\in \ST^H(P\times\{0,1\})$.

\medskip

Theorem \ref{tth} is a corollary of the construction below.  This
construction is almost identical to that in \cite{Jo}, except that we
have removed the use of codimension 1 $\pi$-$\pi$ splitting to gain an
extra dimension $m\geq 6$.  Nonetheless, we describe the proof in
detail here, because we shall need to refer to the details later as we
prove a transverse variation.  A transverse variation of Theorem
\ref{tth} follows from a transverse variation of \ref{construct}.

\begin{prop} \label{construct}
Let $m\geq 6$. Given a connected $m$-dimensional topological manifold $N$ and 
an element $I\in L_0(\Z)$ there exists a degree 1 normal map
$(f_I,b_I):N_I\to N$ from a connected homology manifold $N_I$ such that\,:

\begin{itemize}
\item[(i)] The resolution obstruction of $N_I$ in $L_0(\Z)$ is
$$i(N_I)~=~I \in L_0(\Z)~.$$
\item[(ii)] The composite
$$(N_I)_{TOP}
\xymatrix@C+30pt{\ar[r]^-{(f_{TOP},b_{TOP})}&} 
N_I \xymatrix@C+15pt{\ar[r]^-{(f_I,b_I)}&} N$$
is normally bordant to the identity map.
\item[(iii)] If $M$ is a connected $m$-dimensional homology manifold with
resolution obstruction $i(M)=I \in L_0(\Z)$,
then the composite
$$(M_{TOP})_I\xymatrix@C+15pt{\ar[r]^-{(f_I,b_I)}&}
M_{TOP}\xymatrix@C+30pt{\ar[r]^-{(f_{TOP},b_{TOP})}&} M$$
is normally bordant to the identity map $M\to M$.
\end{itemize}
\end{prop}

\medskip

\noindent{\bf Proof of Theorem \ref{tth}}:
Denote the given function by $\Phi:\TT^H(X)\to \TT(X)\times L_0(\Z)$.
Define a function
$$\begin{array}{l}
\Psi:\TT(X) \times L_0(\Z)\to \TT^H(X)~;\\
((f,b):N \to X,I) \mapsto ((f_I,b_I):N_I \to N 
\xymatrix{\ar[r]^-{(f,b)}&} X)
\end{array}
$$
where $(f_I,b_I):N_I\to N$ is the result of applying proposition \ref{construct}
to the pair $N,I$.

The composition $\Phi\circ\Psi$ is the identity by \ref{construct} 
(ii), and the composition $\Psi\circ \Phi$ is the identity by
\ref{construct} (iii).  \qed 
\medskip

\noindent{\bf Proof of Proposition \ref{construct}}:
\smallskip

\noindent{\bf The construction of $(f_I,b_I):N_I\to N$.}
\smallskip

The construction of $N_I$ is a variation on a construction found in
\cite{BFMW} (Section 7) which is in turn a variation on the main construction
of that paper. In \cite{BFMW} (Section 7) the
construction is performed on a torus, resulting in a homology manifold not homotopy
equivalent to any manifold.  We perform the construction on an arbitrary topological
manifold with (i), (ii) and (iii) above as the result.  
\medskip

Let $\sigma$ denote the element of $H_m(N;\overline{\LL}_{\bullet})$ which
corresponds to the canonical $TOP$ reduction of the Spivak normal fibration
of $N$, with the desired index $I$. Given a sequence $\eta_k>0$
so that $\lim_{k\to \infty}\eta_k=0$.

\smallskip
\noindent{\bf Step I:} Construct a Poincar\'e space $X_0$ with a map
$p_0:X_0\to N$ so that
\begin{itemize}
\item[(i)] $X_0$ is $\eta_0$-Poincar\'e over $N$.
\item[(ii)] $p_0$ is $UV^1$.
\item[(iii)] $p_0$ has controlled surgery obstruction $\sigma \in
H_m(N;\overline{\LL}_{\bullet})$.
\end{itemize}
\medskip

Slice $N$ open along the boundary of a manifold two skeleton, $D$.  So
$N=B\cup_D C$.  We first apply Lemma 4.4 from \cite{BFMW}.  This will
allow us to perform a small homotopy on ${\rm id}_N:N\to N$ to get a new map
$q_0:N \to N$ which restricts to a $UV^1$ map on $B$, $C$ and $D$. 
Because $q_0|_D$ is a $UV^1$ map, the controlled surgery obstruction
group of $D \times [0,1]\stackrel{q_0\text{proj.}}{\to}N$ is
$$L^c\left( \begin{smallmatrix} D \times [0,1] \\
\downarrow \\
N\end{smallmatrix} \right) ~\cong~ H_{m}(N;\overline{\LL}_{\bullet})~\cong~
H_m(N;\LL_{\bullet}) \oplus H_m(N;L_0(\Z))~.$$

Now by Wall realization we construct a normal invariant $\sigma:K \to
D\times [0,1]$ with controlled surgery obstruction $\sigma$,
which is given by a controlled homotopy equivalence $k:D'\longrightarrow D$ on one end and
by the identity on the other.

Gluing $B$ and $C$ back onto $K$ by the identity and by $k$
respectively results in $X_0$ a Poincar\'e complex.  We define a map
$p_0:X_0 \to N$ by $p_0|_{B\cup C}=id$ and $p_0|_K=\sigma$.  By applying
\cite{BFMW} (4.4) we may assume $p_0$ is $UV^1$.  By taking sufficiently fine
control we may assume that $X_0$ is an $\eta_1$-Poincar\'e space over $N$.
\medskip

\begin{figure}
\epsfxsize=4.5in
\center{\leavevmode \epsfbox{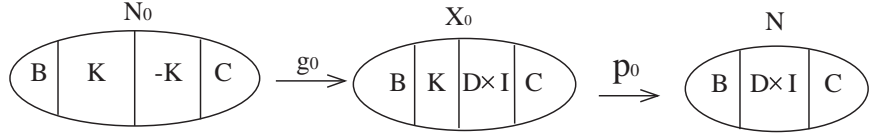}}
\caption{The construction of the spaces $N_0$ 
and $X_0$ and the maps $g_0$ and $p_0$.}
\end{figure} 

The Poincar\'e space $X_0$ has 4-periodic total surgery obstruction
$$\overline{s}(X_0)~=~0\in\overline{\ST}_m(X_0)~,$$
so that $X_0$ is
homotopy equivalent to the desired homology manifold $N_I$
as given by \cite{BFMW} (6.1).  In this variation of the main construction
of \cite{BFMW} the next steps use the Poincar\'e
space $X_0$ and the degree 1 normal map $g_0:N_0\to X_0$ representing
the controlled surgery obstruction $-\sigma$ as described below.
\medskip

Below is a brief summary of the rest of the construction in this case.
It is a limiting process in which the cut and paste
type construction from Step I is performed on finer and finer manifold two
skeleta of manifolds $N_k$.
\medskip

\noindent{\bf Step II:} Construct a Poincar\'e space $X_1$ and a map
$p_1:X_1\to X_0$ so that
\begin{itemize}
\item[\rm (i)] The map $p_1$ is $UV^1$.
\item[\rm (ii)] $X_1$ is an $\eta_1$-Poincar\'e space over $X_0$.
\item[\rm (iii)] The map $p_1:X_1\to X_0$ is an $\epsilon_1$-homotopy
equivalence over $N$.
\item[\rm (iv)] For $W_0$ a regular neighborhood of $X_0$, there exists an
embedding $X_1\to W_0$ and a retraction $r_1:W_0\to X_1$
so that $d(r_0, r_1)<\epsilon_1$.
\end{itemize}
\medskip

Let $N_0$ be given by $B\cup_D K\cup_{D'}-K \cup_{D}C$.  Define a map
$$g_0~:~N_0\to X_0~=~B\cup_D K \cup_{D} D\times [0,1] \cup_D C$$
 by the identity on $B$, $K$ and $C$ and by $-\sigma$ on $-K$.  By
\cite{BFMW} (4.4) we may assume that $g_0$ is $UV^1$.  Let
$N_0=B_1\cup_{D_1}C_1$ be a decomposition of $N_0$ by a finer manifold
two skeleton than that of $N$.  Let $q_0$ denote the map homotopic to
$g_0$ which restricts to a map $UV^1$ on $B_1$, $C_1$ and $D_1$. 
\medskip

Since the map $p_0$ is $UV^1$ it induces an isomorphism
$$(p_0)_*~:~H_m(X_0;\overline{\LL}_{\bullet})~ \cong~
H_m(N;\overline{\LL}_{\bullet})~.$$
Let $\sigma_1:K_1\to D_1\times [0,1]$ denote a Wall realization of the
element of $H_m(X_0;\overline{\LL}_{\bullet})$ which corresponds to
$\sigma$.  Define $X_1'=B_1\cup K_1 \cup C_1$ with a map $f_1:X_1'\to
N_0$ defined as for $p_0$ consider the composition of maps $X_1'\to N_0
\to X_0$ and notice that it has vanishing surgery obstruction and can
therefore be surgered to a small homotopy equivalence over $N$.  (This
type of surgery on a Poincar\'e space is in the tradition of Lowell
Jones \cite{Jn}.) Denote the result of this surgery by $p_1:X_1 \to
X_0$.  We may assume that $p_1$ is $UV^1$.

\medskip

\begin{figure}
\epsfxsize=4.5in
\center{\leavevmode \epsfbox{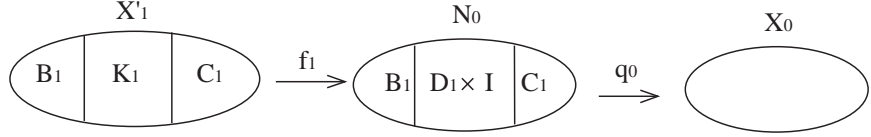}}
\caption{By construction one can surger the composition of maps
$X_1'\stackrel{f_1}{\to}N_0\stackrel{q_0}{\to}X_0$ to a homotopy equivalence
$p_1:X_1\to X_0$.}
\end{figure} 

By choosing a sufficiently well-controlled surgery obstruction, we may
assume that $X_1$ is $\eta_1$-Poincar\'e over $X_0$.  By choosing
$\eta_0$ sufficiently small we may verify conditions (iii) and (iv).
\medskip

\noindent{\bf Step III:} Construct a Poincar\'e space $X_{i+1}$ and a map
$p_{i+1}:X_{i+1}\to X_i$ so that
\begin{itemize}
\item[(i)] $p_{i+1}$ is $UV_1$,
\item[(ii)] $X_{i+1}$ is $\eta_i$-Poincar\'e over $X_i$,
\item[(iii)] $p_{i+1}$ is an $\epsilon_i$ equivalence over $X_{i-1}$,
\item[(iv)] there is an embedding $X_{i+1}\to W_i\subset W_0$ and a retraction
$r_i:W_0\to X_{i+1}$ so that $d(r_i,r_{i+1})<\epsilon_i$.
\end{itemize}

Let $g_i:N_i\to X_i$ be a degree 1 normal map with surgery obstruction
$-\sigma \in H_m(X_{i-1};\overline{\LL}_{\bullet})\cong
H_m(N;\overline{\LL}_{\bullet})$.  By \cite{BFMW} (4.4) we may assume that
$f_i$ is $UV^1$.  Let $N_i=B_i\cup_{D_i}C_i$ be a decomposition
of $N_i$ by a finer manifold two skeleton than that of $N_{i-1}$.  Let
$q_i$ denote the map homotopic to $g_i$ which restricts to a map $UV^1$
on $B_i$, $C_i$ and $D_i$.

\medskip

Since the map $p_i$ is $UV^1$ it induces an isomorphism
$$H_m(X_i;\overline{\LL}_{\bullet})~ \cong~ H_m(X_{i-1};
\overline{\LL}_{\bullet})~\cong~H_m(N;\overline{\LL}_{\bullet})~.$$ 
Let $\sigma_i:N_i\to D_i\times [0,1]$ denote a Wall realization of
the element of $H_m(X_i;\overline{\LL}_{\bullet})$ which corresponds to
$\sigma$.  Define $X_{i+1}'=B_i\cup K_i \cup C_i$ with a map
$f_{i+1}:X_{i+1}'\to N_i$ defined as for $p_0$ consider the composition
of maps $X_{i+1}'\to N_i \to X_i$ and notice that it has vanishing
surgery obstruction and can therefore be surgered to a small homotopy
equivalence over $X_{i-1}$.  Denote the result of this surgery by
$p_{i+1}:X_{i+1} \to X_i$.  We may assume that $p_{i+1}$ is $UV^1$.

\begin{figure}
\epsfxsize=4.5in
\center{\leavevmode \epsfbox{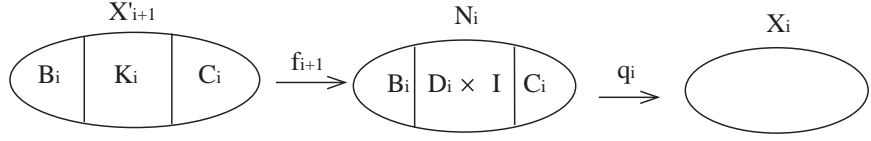}}
\caption{By construction one can surger the composition of maps
$X_{i+1}'\stackrel{f_{i+1}}{\to}N_i\stackrel{q_i}{\to}X_i$ to a homotopy 
equivalence $p_{i+1}:X_{i+1}\to X_i$.}
\end{figure} 
\medskip

By choosing a sufficiently well-controlled surgery obstruction, we may
assume that $X_{i+1}$ is $\eta_{i+1}$-Poincar\'e over $X_i$.  By
choosing $\eta_{i-1}$ sufficiently small we may verify conditions {\rm (iii)}
and {\rm (iv)}.
\medskip

\noindent{\bf Step IV:} Let $N_I=\cap_{i+1}^{\infty}W_i$.  This is the
desired homology manifold and let the map $N_I\to N$ be defined by
$N_I\stackrel{r_k|_{N_I}}{\to} X_0\stackrel{p_0}{\to} N$.
\medskip

$N_I$ is an $ANR$, because the limit of the retractions $r_i$ defines a
retraction $r:W_0\to N_I$.  To see that $N_I$ is a homology manifold,
we first use condition {\rm (iv)} to improve the retractions $r_i$.  Then
this together with the fact that each $X_i$ is an $\eta_i$-Poincar\'e
space over $X_{i-1}$ can be used to show that there is a retraction
$\rho:W_0\to N_I$ so that $\rho|:\partial W_0\to N_I$ is an approximate
fibration, which shows that $N_I$ is a homology manifold.

\medskip

This concludes the construction of $(f_I,b_I):N_I\to N$.
It remains to show that the construction has produced the desired result.

\medskip

\noindent{\bf Proof of (i)}: Consider the controlled surgery
obstruction of $N_k\to X_k\to N_I$ controlled over $N_I$ by the
identity map $N_I\to N_I$.  Since the map $\rho:X_k \to N_I$ can be
assumed to be $UV^1$ this is the same as the controlled surgery obstruction of $g_k:N_k\to
X_k$ where $X_k$ has control map $\rho:X_k\to N_I$.  By the argument
given in \cite{BFMW} the control maps $\rho:X_k\to N_I$ and
$$X_k\stackrel{p_k}{\to}X_{k-1}\stackrel{\rho}{\to} N_I$$
are homotopic by a small homotopy.  Thus the controlled surgery
obstruction of $g_k$ with one control map $\rho$ is the same as that
with control map $p_k$, but by construction this surgery obstruction
was $-\sigma$, i.e.  the resolution obstruction of $N_I$ is $I$.
\qed

\medskip

\noindent {\bf Proof of (ii)}:
Let $f_{TOP}:(N_I)_{TOP}\to N_I$ denote the degree 1 normal map induced
by the $TOP$ reduction of $N_I$ given by the map
$$N_I\stackrel{f_I}{\to} N\stackrel{\nu_N}\to BTOP~.$$
By construction this map is normally bordant to ${\rm id}_N:N\to N$.
\qed
\medskip

\noindent{\bf Proof of (iii)}:
The proof is the same as injectivity of $\Phi$ in Theorem 5.2 of
\cite{Jo}, except that we now restrict our attention to showing that
the composition $(f_{TOP})_I\circ f_{TOP}:(M_{TOP})_I\to M_{TOP}\to M$
is homology manifold normal bordant to the identity map ${\rm id}_M:M\to M$. 
Since any map $f_I:N_I\to N$ is an isomorphism on fundamental groups,
we can avoid Lemma 5.4 of \cite{Jo} and its requirement that $m\geq 7$. 
The map $f_I$ is an isomorphism on fundamental groups, because it
splits as a homotopy equivalence $h:N_I\to X_0$ and the original map
$p_0:X_0\to N$.  The fact that $p_0$ is an isomorphism on fundamental
groups follows from the fact that the normal invariant $\sigma:K\to
D\times [0,1]$ is a Wall realization and therefore an isomorphism on 
fundamental groups.

A homology manifold normal bordism between $(f_{TOP})_I\circ
f_{TOP}:(M_{TOP})_I\to M_{TOP}\to M$ and ${\rm id}_M:M\to M$ is constructed
as follows. First, a Poincar\'e bordism $k:W\to M$ is constructed
between the two maps.  Then using the fact that we may assume that
$(W,(M_{TOP})_I)$ is a $(\pi,\pi)$ pair we put a homology manifold
structure on $W$ relative to $M$.  See \cite{Jo} (5.2) for details.  
\qed

\medskip

For any topological block bundle $\nu:X \to B\widetilde{TOP}(k)$
define the homology manifold bordism groups $\Omega^H_m(X,\nu)$
of normal maps $(f,b):M \to X$ from $m$-dimensional
homology manifolds, with $b:\nu_M \to \nu$.

\begin{cor} \label{jbord}
For $m\geq 6$ the $m$-dimensional homology
manifold bordism groups of $(X,\nu)$ are such that
$$\Omega^H_m(X,\nu)~\cong~\Omega^{TOP}_m(X,\nu)[L_0(\Z)]~.$$
\end{cor}

\begin{proof} Use the construction of Proposition \ref{construct} to define
inverse isomorphisms
$$ \begin{array}{l}
\psi~:~\Omega^H_m(X,\nu) \to \Omega^{TOP}_m(X,\nu)[L_0(\Z)]~;\\
\hskip75pt
((f,b):M \to X) \mapsto (M_{TOP} \to M \stackrel{(f,b)}{\to} X)[i(M)]~,\\[2ex]
\psi^{-1}~:~\Omega^{TOP}_m(X,\nu)[L_0(\Z)]\to \Omega^H_m(X,\nu)~;\\
\hskip75pt
((g,c):N \to X)[I] \mapsto (N_I \to N \stackrel{(g,c)}{\to}  X)~.
\end{array}$$
\end{proof}

Thus
$$\Omega^H_m(X,\nu)~=~\Omega^{TOP}_m(X,\nu)[L_0(\Z)]~=~\pi_{m+k}(T(\nu))[L_0(\Z)]~,$$
with $T(\nu)$ the Thom space. In particular, for any space $K$ and
$$\nu~=~{\rm proj.}~:~X~=~K\times BTOP \to BTOP$$
we have
$$\Omega^H_m(K)~=~\Omega^{TOP}_m(K)[L_0(\Z)]~=~
H_m(K;\Omega^{TOP}(\{{\rm pt.}\}))[L_0(\Z)]~.$$
If $X$ is an $m$-dimensional Poincar\'e space
with $TOP$ reducible Spivak normal fibration then for
each $TOP$ reduction $\nu:X \to B\widetilde{TOP}(k)$
$$\Omega^H_m(X,\nu)~=~\pi_{m+k}(T(\nu))[L_0(\Z)]~=~[X,G][L_0(\Z)]~.$$
For connected $X$ this is another way to see that
$$\TT^H(X)~=~\TT^{TOP}(X) \times L_0(\Z)~=~[X,G/TOP]\times L_0(\Z)$$
as in Theorem \ref{eleven}.

\section{Homology submanifold transversality up to bordism}

Fix a space $P$ with a codimension $q$ bundle subspace $(Q,R,\xi)$ as
in the Introduction.
\medskip

We now investigate the transversality to $Q \subset P$
of a map $f:M^m \to P$ from an $m$-dimensional homology manifold.
In the first instance, we show that if $m-q\geq 7$ then
$f$ is bordant to a transverse map.

\begin{definition} {\normalfont
The {\it $Q$-transverse homology manifold bordism group}
$\Omega^{H,Q-\transverse}_m(P)$ is the abelian group of bordism classes
of maps $M^m \to P$ from $m$-dimensional homology manifolds which are
transverse to $Q \subset P$.}
\end{definition}

There are evident forgetful maps
$\Omega^{H,Q-\transverse}_m(P)\to \Omega^H_m(P)$.

\begin{theorem} \label{eight}
{\rm (i)} The $Q$-transverse homology manifold bordism
groups fit into an exact sequence
$$\dots \to\Omega^H_m(R) \to \Omega^{H,Q-\transverse}_m(P) \to
\Omega^H_{m-q}(Q) \to \Omega^H_{m-1}(R)\to \dots$$
with
$$\begin{array}{l}
\Omega^H_m(R) \to \Omega^{H,Q-\transverse}_m(P)~;~(M,f:M \to R) \mapsto
(M,M \stackrel{f}{\to} R \to P)~,\vspace*{2mm}\\
\Omega^{H,Q-\transverse}_m(P)\to \Omega^H_{m-q}(Q) ~;~
(M,g:M \to P) \mapsto (g^{-1}(Q),g\vert:g^{-1}(Q) \to Q)~,\vspace*{2mm}\\
\Omega^H_{m-q}(Q) \to \Omega^H_{m-1}(R)~;~
(N,h:N \to Q) \mapsto (S(h^*\xi),S(h^*\xi) \to S(\xi) \to R)~.
\end{array}$$
{\rm (ii)} For $m-q \geq 6$
$$\Omega^{H,Q-\transverse}_m(P)~=~\Omega^H_m(P)~=~
\Omega^{TOP}_m(P)[L_0(\Z)]~.$$
In particular, the forgetful maps
$\Omega^{H,Q-\transverse}_m(P)\to \Omega^H_m(P)$
are isomorphisms, and every map $M^m \to P$ is bordant to a $Q$-transverse map.
\end{theorem}
\begin{proof} {\rm (i)} This is a formality.\\
{\rm (ii)} Define the $Q$-transverse topological bordism groups
$\Omega^{TOP,Q-\transverse}_m(P)$ by analogy with
$\Omega^{H,Q-\transverse}_m(P)$, for which there is an exact sequence
$$\dots \to\Omega^{TOP}_m(R) \to \Omega^{TOP,Q-\transverse}_m(P) \to
\Omega^{TOP}_{m-q}(Q) \to \Omega^{TOP}_{m-1}(R)\to \dots~.$$
The forgetful maps $\Omega^{TOP,Q-\transverse}_m(P) \to \Omega^{TOP}_m(P)$
are isomorphisms, by topological transversality. Applying the
5-lemma to the map of exact sequences
$$\xymatrix@C1pc{\dots \ar[r] & \Omega_m^H(R) \ar[d]^{\cong} \ar[r] &
\Omega^{H,Q-\pitchfork}_m(P)
\ar[d]\ar[r] & \Omega^H_{m-q}(Q) \ar[r] \ar[d]^{\cong} & \dots\\
\dots \ar[r] & \Omega_m^{TOP}(R)[L_0(\Z)] \ar[r] &
\Omega^{TOP,Q-\pitchfork}_m(P)[L_0(\Z)]
\ar[r] & \Omega^{TOP}_{m-q}(Q)[L_0(\Z)]  \ar[r] & \dots
} $$
we have that the morphisms
$$\Omega^{H,Q-\transverse}_m(P) \to \Omega^{TOP,Q-\transverse}_m(P)[L_0(\Z)]~;~
(M,f:M \to P) \mapsto (M_{TOP},f_{TOP})[i(M)]$$
are isomorphisms for $m-q \geq 6$, exactly as in the case $Q=\emptyset$.
\end{proof}

We wish to improve this homology manifold transversality result from
``up to bordism" to ``up to normal cobordism", using surgery theory.
\medskip

We shall need the following variation of \ref{construct}\,:
\begin{prop}\label{qconstruct}
Let $m-q\geq 7$. Given a connected $m$-dimensional topological manifold $N$,
a $Q$-transverse map $g:N\to P$, and an
element $I\in L_0(\Z)$ there exists a degree 1 normal map
$(f_I,b_I):N_I\to N$ from a connected homology manifold $N_I$ such that\,:
\begin{itemize}
\item[(i)] The resolution obstruction of $N_I$ in $L_0(\Z)$ is
$$i(N_I)~=~I \in L_0(\Z)~.$$
\item[(ii)] The composite
$$(N_I)_{TOP} \xymatrix@C+30pt{\ar[r]^-{(f_{TOP},b_{TOP})}&} 
N_I\xymatrix@C+10pt{\ar[r]^-{(f_I,b_I)}&} N$$
is normally bordant to the identity map by a normal bordism.
(By a $Q$-transverse normal bordism, since we are in the topological category.)
\item[(iii)] If $M$ is a connected $m$-dimensional homology manifold
with resolution obstruction $i(M)=I \in L_0(\Z)$, then the composite
$$(M_{TOP})_I \xymatrix@C+10pt{\ar[r]^-{(f_I,b_I)}&}
M_{TOP} \xymatrix@C+30pt{\ar[r]^-{(f_{TOP},b_{TOP})}&} M$$ 
is normally bordant to the identity map $M\to M$.
\item[(iv)] The map $(f_I,b_I)$ is $Q$-transverse.
\end{itemize}
\end{prop}

\begin{proof}
This result follows from the following two lemmas, whose proofs we defer.

\begin{lemma}\label{uconstr}
Given an $m$-dimensional topological manifold $M$ for $m\geq 7$, such that $M$ is the union
of two manifolds
along a boundary component, $M=M_+\cup_{M_0}M_-$.
We may perform the construction of \ref{construct} so that the result is a
homology manifold $M_I=(M_+)_I\cup_{(M_0)_I} (M_-)_I$ so that $M_I$ and
$(M_i)_I$ for $i=0,+,-$
all satisfy the conclusion of \ref{construct}.
\end{lemma}

\begin{lemma}\label{bundleconstr} Given an $(m-q)$-dimensional manifold $X$ for $m-q\geq 7$
and a $D^q$-bundle
$\xi$ over $X$, whose total space is $E(\xi)$ we may construct a
homology manifold $E_I$ with a map $$f_I:E_I\to E(\xi)$$ which satisfies
the properties of \ref{qconstruct}.  In particular $f_I^{-1}(X)=X_I$ is
a homology manifold such that $E_I$ is the total space of the bundle
$f_I^*\xi$ and $$f_I|:X_I\to X$$ also satisfies the properties of
\ref{construct}.  \end{lemma}

Denote $g^{-1}(Q)$ by $N_Q$ and $g^{-1}(R)$ by $N_R$, so that
$(N_Q,N_R,g^*\xi)$ is a codimension $q$ bundle subspace of $N$ with
$$N=E(g^*\xi)\cup N_R.$$  By applying Lemma \ref{bundleconstr} to
$E(g^*\xi)$, we get $$(f_I)_E:E_I\to E(g^*\xi)$$ satisfying the above
conditions.  Applying \ref{construct} to $N_R$ results in
$$(f_I)_R:R_I\to N_R.$$  By Lemma \ref{uconstr} we may apply the
construction to $N_R$ and $E(g^*\xi)$ simultaneously so that the
resulting homology manifolds $R_I$ and $E_I$ and maps agree on their
boundaries.  \end{proof}

\medskip

\noindent{\bf Proof of Lemma \ref{uconstr}:}
Take a manifold two-skeleton of $M$ which is the union of manifold two
skeletons for $M_+$ and $M_-$ along a manifold two skeleton for $M_0$.
Denote this two skeleton and its boundary by
$$B=B_+\cup_{B_0}B_- \text{ and }D=D_+\cup_{D_0}D_-.$$  We shall perform the
construction of \ref{construct} simultaneously on $M_+$ and $M_-$. 
What was one controlled surgery obstruction in the original
construction
$$\sigma=(0,I)\in H_m(M;\overline{\LL}_{\bullet})$$ 
is now two controlled surgery obstructions
$$\sigma_+ ~=~ (0,I)\in H_m(M_+;\overline{\LL}_{\bullet})~~,~~
\sigma_-~=~(0,I)\in H_m(M_-;\overline{\LL}_{\bullet})~.$$
\medskip

In Step I, we apply \cite{BFMW} (4.4) to get a map $q_0:M\to M$ which
restricts to $UV^1$ maps $D_i\to M_i$,~$i=0,+,-$ and is itself $UV^1$. 
Thus the controlled surgery obstruction groups are given by
$$L^c\left( \begin{smallmatrix} D_i \times [0,1] \\
\downarrow \\
M_\pm\end{smallmatrix} \right) ~\cong~ H_{m}(M_\pm;\overline{\LL}_{\bullet})~\cong~
H_m(M_\pm;\LL_{\bullet}) \oplus H_m(M_\pm;L_0(\Z))~.$$
Since $M_\pm$ are manifolds with boundary, $m\geq 7$ is required.
Take a Wall realization $K_\pm$ of $(0,I)$.  Glue this into $M_\pm$ to
create the first Poincar\'e space of the construction, which comes with
a map $p_0:X_0\to M$ and a Poincar\'e decomposition
$X_0=(X_0)_+\cup(X_0)_-$ which is respected by $p_0$.
\medskip

We must exercise some care to get the Wall realizations to agree on
their boundaries.
\begin{figure}
\epsfysize=1in
\center{\leavevmode \epsfbox{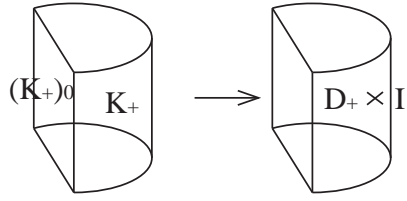}}
\caption{The Wall realization $K_+$.}
\end{figure} 

Denote $\partial K_\pm=(K_\pm)_0$.  The images of $(K_\pm)_0\to D_0\times [0,1]$
agree in
$$L^c\left( \begin{smallmatrix} D_0 \times [0,1] \\ \downarrow
\\ M_0\end{smallmatrix} \right) ~\cong~
H_{m-1}(M_0;\overline{\LL}_{\bullet})~.$$
Each is again $(0,I)$, the image of $(0,I)$ under the boundary maps
$$H_m(M_\pm;\overline{\LL}_{\bullet})\to H_{m-1}(M_0;\overline{\LL}_{\bullet})~.$$
Since the surgery obstructions of $\partial K_+=(K_+)_0$ and $\partial K_-=(K_-)_0$
agree we may glue them together along their common boundary
$D_0$ to get a normal invariant with vanishing controlled
surgery obstruction, i.e.  we may surger the map $$\partial K_+\cup\partial
K_-\to D_0\times[0,1]$$ to a controlled homotopy equivalence. Let $W$ denote
the trace of this surgery.
Denote $K_-'=K_-\cup_{(K_-)_0} W$, and extend the map $K_-\to D_-\times [0,1]$ in the
obvious way.
It is a Wall realization of $(0,I)$, which agrees with $K_+$ on its
boundary.  \medskip

\begin{figure}
\epsfxsize=3in
\center{\leavevmode \epsfbox{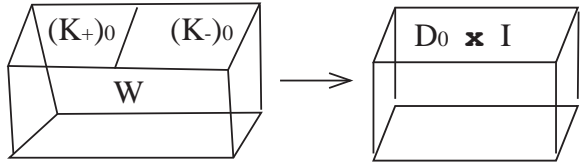}}
\caption{The trace of surgery on $(K_+)_0\cup (K_-)_0\to D_0\times[0,1]$ 
can be used to ``match up'' the boundaries of the Wall realizations.}
\end{figure}

If we perform the construction with two Wall realizations which agree
on their boundaries, then we may glue them together to get
$X_0=(X_0)_+\cup(X_0)_-$.  Furthermore since the boundary of the Wall
realizations is itself a Wall realization, the intersection
$$(X_0)_+\cap(X_0)_-$$ is itself the first stage of the given
construction on $M_0$.  Similarly the union $X_0$, which is the result
of gluing in the union of the Wall realizations is the first stage of
the construction on $M$.  \medskip

To preserve this decomposition throughout the construction requires
only that we repeat the above type of construction in later steps.
\qed \medskip

\noindent {\bf Proof of Lemma \ref{bundleconstr}:}
Let $B_X$ denote a manifold two skeleton of $X$.  Since $E(\xi)$ and
$X$ are of the same homotopy type we could easily construct a manifold
two skeleton of $E(\xi)$ by thickening a manifold two skeleton of $X$,
i.e.  $E(c^*\xi)$ where $c:B_X\to X$.  Recall that the construction
requires a fine manifold two skeleton.  Let $E_X$ denote the manifold
two skeleton of an $\epsilon$ neighborhood of $X$ in $E(\xi)$.  We may
extend this to a fine manifold two skeleton of $E(\xi)$ by 
$$B~=~E_X\cup B_S$$ 
where $B_S$ is a manifold two skeleton of the complement of the
given $\epsilon$ neighborhood of $X$.  
\medskip

Here $E_X\cap B_S$ is the manifold two skeleton of $S(\xi)$ for the
given $\epsilon$-sphere bundle.  We shall now perform the construction
of \ref{uconstr} on this decomposition of $E(\xi)$, but with the added
requirement that the result of the construction on the small
neighborhood of $X$ is itself a bundle over the desired $X_I$.  This is
done by preserving the bundle throughout the construction.  Take a Wall
realization of 
$$(0,I)\in H_{m-q}(X;\overline{\LL}_{\bullet})$$ 
denoted by
$$\sigma_X:N_X\to D_X\times [0,1]~.$$ 
Let $\sigma_E$ denote the bundle map on the pullback 
$$N_E~=~\sigma_X^*\xi'$$ 
where $\xi'$ is the bundle induced on $D_X\times[0,1]$ by $\xi$. 
The result is
$$\sigma_E:N_E\to E(\xi|D_X\times [0,1])~,$$ 
the image of $\sigma_X$ under the transfer map.  The corresponding surgery
obstruction is thus the image of $(0,I)$ under the transfer map
$$H_{m-q}(X;\overline{\LL}_{\bullet})\to
H_m(E(\xi);\overline{\LL}_{\bullet})~.$$ 
One need only check that this map takes $(0,I)$ to $(0,I)$ to see that
we may preserve the bundle throughout the construction.
\qed 

\medskip

Given $X$ a simple $m$-dimensional Poincar\'e space with a map $g:X \to
P$ to a space $P$ with a codimension $q$ bundle subspace $Q \subset P$.  Let
$\TT^{H,Q-\transverse}(X)$ be the set of bordism classes of normal maps
$(f,b):M \to X$ from $m$-dimensional homology manifolds such that $gf:M
\to P$ is transverse to $Q \subset P$.

\begin{theorem}\label{bordtransverse} Given $X$, $g:X\to P$ and $Q\subset P$
as above such that $m-q\geq 7$.
\\
{\rm (i)} The forgetful function $\TT^{H,Q-\transverse}(X) \to \TT^H(X)$
is a bijection.\\
{\rm (ii)} Every map $f:M \to P$ from an $m$-dimensional homology manifold
is bordant to a map transverse to $Q \subset P$.\end{theorem}

\begin{proof} (i) The functions
$$\begin{array}{l}
\TT^H(X) \to \TT^{TOP}(X)\times L_0(\Z)\,;\,((f,b):M \to X) \mapsto
((f_{TOP},b_{TOP}),i(M))\,,
\vspace*{2mm}\\
\TT^{H,Q-\transverse}(X) \to \TT^{TOP,Q-\transverse}(X)\times L_0(\Z)\,;\,
((f,b):M \to X) \mapsto ((f_{TOP},b_{TOP}),i(M))
\end{array}$$
are bijections by Theorem \ref{tth} and its $Q$-transverse variation
\ref{qconstruct}.
The forgetful function
$$\TT^{TOP,Q-\transverse}(X) \to \TT^{TOP}(X)$$
is a bijection by topological transversality, so that
$$\TT^{H,Q-\transverse}(X)~=~\TT^{TOP,Q-\transverse}(X)\times L_0(\Z)
~=~\TT^{TOP}(X)\times L_0(\Z)~=~\TT^H(X)~.$$
(ii) Unfortunately (i) does not apply to an arbitrary map.  We may get
around this by factoring any map $f:M\to P$ through a homotopy
equivalence $\widehat{f}:M\to \overline{P}$.  Any map $f$ is homotopic to
$\widehat{f}\circ\overline{f}$ such that $\widehat{f}:\overline{P}\to P$ is a Serre
fibration and $\overline{f}$ is a homotopy equivalence.  Now by (i),
$\overline{f}$ is normally cobordant to a transverse map, and hence $f$ is
bordant to a transverse map.  
\end{proof}

Given $X$ a simple $m$-dimensional Poincar\'e duality space with a map
$g:X \to P$, as above, we may also study the {\it $Q$-transverse
homology manifold structure set}.  Denoted by,
$\ST^{H,Q-\transverse}(X)$, this is the set of equivalence classes of
pairs $(M,h)$ with $M$ an $m$-dimensional homology manifold and $h:M
\to X$ a simple homotopy equivalence such that $gh:M \to P$ is
transverse to $Q\subset P$, with $(M_1,h_1) \allowbreak \simeq
(M_2,h_2)$ if there exists an $s$-cobordism $(W;M_1,M_2)$ with a simple
homotopy equivalence of the type
$$(f;h_1,h_2)~:~(W;M_1,M_2) \to X \times ([0,1];\{0\},\{1\})$$
such that the composite
$$W \stackrel{f}{\to} X \times [0,1] \stackrel{proj.}{\to} X
\stackrel{g}{\to} P$$
is transverse to $Q\subset P$.

\begin{remark}
The isomorphism 
$$\ST^{H,Q-\transverse}(X)~\cong~\ST^H(X)$$ 
will follow from $s$-transversality for homology submanifolds, Theorem
\ref{homtra}, just as Theorem \ref{bordtransverse} followed from
\ref{qconstruct}.  
\end{remark}

Before dealing with transversality up to $s$-cobordism, we turn our
attention to some splitting theorems.  These theorems follow directly
from transversality up to normal bordism and are useful in the proof of
transversality up to $s$-cobordism.

\section{Codimension $q$ splitting of homology manifolds}

Fix a space $P$ with a codimension $q$ bundle subspace $(Q,R,\xi)$,
as in Section 2.
\medskip

Wall \cite{Wa} (Chapter 11) defined the codimension $q$ splitting obstruction
groups $LS_*(\Phi)$ to fit into an exact sequence
$$\begin{array}{rl}
\dots \to L_{m+1}(\Z[\pi_1(R)]\!\!\!&\to\Z[\pi_1(P)])\to LS_{m-q}(\Phi)\\
&\to L_{m-q}(\Z[\pi_1(Q)])
\stackrel{\xi^!}{\to} L_m(\Z[\pi_1(R)]\to\Z[\pi_1(P)])\to \dots
\end{array}$$
with $\xi^!$ the transfer maps induced by $\xi$.

\medskip

From now on, we assume that $P$ is an
$m$-dimensional Poincar\'e space and that $Q \subset P$ is an
$(m-q)$-dimensional Poincar\'e subspace, with $(R,S(\xi))$ an
$m$-dimen-sional Poincar\'e pair.

\begin{definition} {\normalfont
{\rm (i)} A simple homotopy equivalence $f:M \to P$ from an $m$-dimensional homology
manifold
{\it splits along} $Q \subset P$  if $f$ is $s$-cobordant to a simple homotopy
equivalence
(also denoted by $f$) which is transverse to $Q \subset P$
$$f~=~E(g)\cup_{S(g)}h~:~M~=~E(g^*\xi) \cup_{S(g^*\xi)}Z \to P~=~E(\xi)\cup_{S(\xi)}R$$
and such that the restrictions
$$g~=~f\vert~:~N ~=~f^{-1}(Q) \to Q~,~h~=~f\vert~:~Z ~=~f^{-1}(R) \to R$$
are simple homotopy equivalences.\\
{\rm (ii)} The {\it split structure set} $\ST^H(P,Q,\xi)$
is the set of homology manifold structures on $P$
which split along $Q \subset P$.}
\end{definition}

For any simple homotopy equivalence
$f:M \to P$ from an $m$-dimensional topological manifold $M$
there is
defined a codimension $q$ splitting obstruction
$$s_Q(f) \in LS_{m-q}(\Phi)$$
such that $s_Q(f)=0$  if (and for $m-q \geq 5$ only if) $f$ splits in the
topological manifold category. The image of $s_Q(f)$ in
$L_{m-q}(\Z[\pi_1(Q)])$ is the surgery obstruction $\sigma_*(g)$
of the normal map $g=f\vert:N=f^{-1}(Q) \to Q$ obtained by topological
transversality. See Wall \cite{Wa} (Chapter 11) 
and Ranicki \cite{R2} (pp.\,572--577) for details.

\medskip

\begin{theorem} \label{qsplit} Let $m-q \geq 7$.\\
{\rm (i)} A simple homotopy equivalence $f:M \to P$ from an $m$-dimensional homology
manifold with a codimension $q$ bundle subspace $(Q, R, \xi)$ splits along $Q \subset P$
if and only if an obstruction $s_Q(f) \in LS_{m-q}(\Phi)$ vanishes.\\
{\rm (ii)} The split homology manifold
structure set $\ST^H(P,Q,\xi)$ fits into an exact sequence
$$\dots \to LS_{m-q+1}(\Phi)\to \ST^H(P,Q,\xi) \to \ST^H(P) \to
LS_{m-q}(\Phi)~.$$
\end{theorem}
\begin{proof}
Define $LS^H_{m-q}(\Phi)$ to be the group of obstructions in the exact sequence
$$\ST^H(P,Q,\xi)\to NI^{H,Q-\transverse}(X)\to LS^H_{m-q}(\Phi)~.$$

Consider the homology manifold normal invariant given by $f:M\to P$, by
\ref{bordtransverse} $f$ is normally bordant, say via $W\to P\times I$
to a map $g:M'\to P$ such that $g$ is transverse to $Q$.  In particular
$g^{-1}(Q)=Y$ is a homology manifold with a normal neighborhood
$N(Y)=g^*(\xi)$.

Thus, this transverse normal invariant defines a splitting obstruction
which lives in $LS^H_{m-q}(\Phi)$.  Since we do not a priori have an
understanding of $LS^H_{m-q}(\Phi)$ we must study it by comparing to
the obstruction groups $L_{m+1}(\pi_1(R)\to \pi_1(P)), L_{m-q}(Q)$ and
the surgery exact sequences for $\ST^H(P\times[0,1],R\times 0)$ and
$\ST^H(Q)$.

There is clearly a commutative diagram with vertical maps given by
restriction as follows:
\begin{equation*}
\begin{CD}
\ST^H(P\times I,R\times 0)@>>> NI^H(P\times I,R\times 0)@>>>
L_{m+1}(R\to P)\\
@VVV @VVV @VVV\\
\ST^H(P,Q,\xi)@>>>NI^{H,Q-\transverse}(P)@>>>LS^H_{m-q}(\Phi)\\
@VVV @VVV @VVV\\
\ST^H(Q)@>>> NI^H(Q)@>>>L_{m-q}(Q)\\
\end{CD}
\end{equation*}

The splitting obstruction $\sigma(f)\in LS^H_{m-q}(\Phi)$ can be
understood as a two stage obstruction as follows.  First the normal
invariant $g|:Y\to Q$ defines an obstruction $\sigma_Q(f)\in
L_{m-q}(Q)$.  If this obstruction vanishes, then by the surgery exact
sequence of \cite{BFMW} $g|:Y\to Q$ is normally cobordant to a simple
homotopy equivalence.  Let $V\to Q\times [0,1]$ denote this normal
cobordism, and $N(V)\to E(\xi)\times [0,1]$ denote the corresponding
pullback of $\xi$.  We can define a homology manifold normal invariant
of $P\times [0,1] \text{ rel} P\times 0, E(\xi)\times 1$ by gluing
$N(V)$ to $W$.  There is now defined an obstruction $\sigma_R(f)\in
L_{m+1}(\pi_1(R)\to \pi_1(P))$ to this normal invariant being
equivalent to a simple homotopy equivalence of $P\times [0,1]$, i.e. 
an s-cobordism from $f:M\to P$ to some map $h:N\to P$.

\begin{figure}
\epsfysize=.75in
\center{\leavevmode \epsfbox{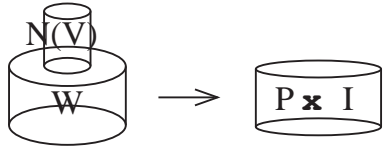}}
\caption{A homology manifold normal invariant of $(P\times[0,1], R\times 0)$.}
\end{figure} 

If $\sigma(f)\in LS_{m-q}(\Phi)$ vanishes, then both $\sigma_Q(f)$ and
$\sigma_R(f)$ are defined and vanish, so that $f$ is $s$-cobordant to a
split map.  Conversely if $f$ is $s$-cobordant to a split map, then
both of $\sigma_Q(f)$ and $\sigma_R(f)$ are defined and vanish.
\end{proof}

\begin{example} {\normalfont
{\rm (i)} If $q=1$, $\xi$ is trivial, and $R=R_1 \sqcup R_2$ is
disconnected with $\pi_1(Q)\cong \pi_1(R_1)$,
then $LS_*(\Phi)=0$, there is no
obstruction in Theorem \ref{qsplit}, and
$$\ST^H(P)~=~\ST^H(P,Q,\xi)~.$$
This is the homology manifold $\pi$-$\pi$ codimension 1
splitting theorem already obtained by Johnston \cite{Jo}.\\
{\rm (ii)} If $q \geq 3$ the codimension $q$ splitting obstruction is just
the ordinary surgery obstruction
$$s_Q(f)~=~\sigma_*(g) \in LS_{m-q}(\Phi)~=~L_{m-q}(\Z[\pi_1(Q)])$$
of the restriction $g=f\vert : f^{-1}(Q) \to Q$,
with an exact sequence
$$\dots \to L_{m-q+1}(\Z[\pi_1(Q)])\to \ST^H(P,Q,\xi) \to \ST^H(P) \to
L_{m-q}(\Z[\pi_1(Q)])~.$$
This is the homology manifold Browder splitting theorem already
obtained by Johnston \cite{Jo}.}
\end{example}
\qed

\section {Homology submanifold transversality up to $s$-cobordism}

We proceed to prove homology manifold transversality up to $s$-cobordism, using the
above results.

\begin{theorem}\label{homtra}
Let $f:M \to P=E(\xi)\cup_{S(\xi)}R$ be a map from an $m$-dimensional
homology manifold $M$ to a space $P$ with a codimension $q$ bundle
subspace $(Q,R,\xi)$.  If $m-q\geq 7$ then $f$ is $s$-transverse to 
$Q \subset P$.
\end{theorem}
\begin{proof} For $Q\subset P$ of codimension $q \geq 3$ this was proved
in \cite{Jo}.

\medskip
We prove the theorem here for codimension $q=1,2$. First we may assume
that the map $f:M\to P$ is a homotopy equivalence by factoring the
original map through a Serre fibration. This results in
$$M\stackrel{\widehat{f}}{\longrightarrow}\overline{P}
\stackrel{\overline{f}}{\longrightarrow}P~,$$ 
homotopic to the original $f$ so that $\widehat{f}$ is a homotopy
equivalence, $\overline{f}$ is a Serre fibration, and $\overline{P}$ has a
codimension $q$ subset $\overline{Q}$ so that the normal bundle of $\overline{Q}$
in $\overline{P}$ is the pullback of the normal bundle of $Q$ in $P$.  To
achieve $Q$-transversality for $f$ clearly it would suffice to achieve
$\overline{Q}$ transversality for $\widehat{f}$.  Thus we may assume that
$f:M\to P$ is a homotopy equivalence.  By theorem \ref{bordtransverse}
(ii) we have that $f$ is bordant to a $Q$-transverse map.  From here
the proof proceeds in two steps, given by the following lemma.

\begin{lemma}\label{lem1}
If $f:M\to P$ is a homotopy equivalence as above and 
$$F~:~W\to P\times [0,1]$$ 
is a homology manifold bordism from $f$ to $g:M'\to P$, then \\
{\rm (i)} The map $f$ is homotopy equivalent to a map $\overline{g}:M\to P$,
which factors through a homotopy equivalence to a Poincar\'e space $X$,
$h:M\to X$ with $g'':X\to P$ such that $\overline{g}=g''\circ h$, and $g''$
is Poincar\'e transverse to $Q$, i.e.  the inverse image of $Q$ is a
Poincar\'e space with normal bundle the pullback of the normal bundle
of $Q$.\\
{\rm (ii)} The map $\overline{g}$ is $s$-cobordant to a $Q$-transverse map.
\end{lemma}
\end{proof}
{\bf Proof of Lemma \ref{lem1} (i):} The key idea of this proof
is to use a patch space structure on $W$ to achieve the desired result. 
Because $W$ is a homology manifold, it has a patch space structure with
only two patches.  Let 
$$H~:~W'\to W$$ 
be such a structure, where
$$W'~=~W_+\cup_{W_0}W_-$$ 
gives $W$ as a union of manifolds glued along a homotopy equivalence 
$$h_W~:~W_0\to W_0~.$$ 
Furthermore we may assume that
$$\partial W'~=~X\amalg X'$$ 
is such that
$$H|~:~\partial W'\to M\amalg M'$$ 
restricts to two patch space structures 
$$X~=~X_+\cup X_-$$ 
and
$$X'~=~X'_+\cup X'_-$$ 
with gluing maps 
$$h_X~:~X_0\to X_0$$ 
and
$$h_{X'}~:~X'_0\to X'_0$$ 
on each of $M$ and $M'$.  
\medskip

We construct a patch space structure for $M'$
as follows: First construct a patch space structure
$$H|~:~Q''\to Q'$$ 
for the homology manifold $Q'=g^{-1}(Q)$.  Thicken this patch space
structure to a patch space structure for $\nu(Q')$ the pullback of the
normal bundle of $Q$.  Then construct a patch space structure for
$$M'\setminus \nu(Q')$$ 
relative to the structure for $\partial\nu(Q')$.  Finally construct a
patch space structure for $W$ relative to this structure for $M'$.  In
this manner a patch space structure for $W$ is constructed so that the
map
$$g\circ H|~:~X'\to P\times \{1\}$$ 
is Poincar\'e transverse to $Q$.
\medskip

Consider the map 
$$G~=~F\circ H~:~W'\to W\to P\times [0,1]~.$$  
We shall show that $G$ is homotopic rel $X'$ to a map which is
Poincar\'e transverse to $Q$.  First consider $G|W_+$. We may use
manifold transversality to make this map transverse to 
$$Q\times [0,1]\subset P\times [0,1]~,$$ 
which inherits a codimension-$q$ structure from $Q\subset P$. If the homotopy 
equivalence
$$h_W~:~W_0\to W_0$$ 
splits along 
$$Q_0~=~(G|_{W+})^{-1}(Q\times [0,1])\cap W_0~,$$ 
then $G|W_-$ is $Q\times [0,1]$-transverse along its boundary $W_0$ and
we may use manifold transversality to homotope $G|W_-$ rel $W_0$ to a
$Q$-transverse map.  Unfortunately there is no guarantee that the
homotopy equivalence $h_W$ splits along $Q_0$.  There is a priori a
splitting obstruction $\sigma \in LS_{m-q}(\Phi)$ where $\Phi$ is
defined by the pair $Q_0\subset W_0$.
\medskip

\begin{figure}
\epsfysize=.75in
\center{\leavevmode \epsfbox{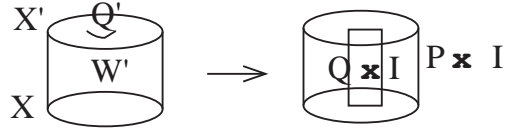}}
\caption{We would like to extend Poincar\'e transversality for the map $G$ 
from $X'$ to all of $W'$.}
\end{figure} 

We shall show that by changing to a different patch space structure
$\overline{H}:\overline{W}\to W$, we may assume that this obstruction vanishes. 
First we need to change the patch space structure slightly so that
$LS_{m-q}(\Phi)$ acts on the set of possible patch space structures. 
Our action will be by $LS_{m-q}(\Psi)$ where $\Psi$ is defined by the
pair $Q''\subset X'$.
\medskip

The first step is to show that we may assume 
$$LS_{m-q}(\Phi)~\cong~LS_{m-q}(\Psi)~.$$  
This is achieved by performing Poincar\'e surgery
to fix the fundamental groups of $Q_0$ and $W_0$, but these are the
boundaries of the manifold two skeletons of 
$$Q'~=~(G|_{W_+})^{-1}(Q \times [0,1])$$ 
and $W'$.  Therefore the fundamental groups of
$Q'\subset W'$ agree with $\Phi$ and we may view the problem of
changing $\Phi$ to agree with $\Psi$ as a problem about $Q'\subset
W'$.  This is just fundamental group Poincar\'e surgery on 
$$G~:~W'\to P\times [0,1]~,$$ 
or rather manifold surgery on the manifold two skeleton of $W'$
$$G|~:~W_+\to W_+~.$$ 
\medskip

Now we can assume that both $Q''\subset X'$ and $Q_0\subset W_0$ define
the same group $LS_{m-q}(\Phi)$.  This group acts on the patch space
structures of $X'$ by acting on $h_{X'}$ as follows: $L_{m-q}(Q'')$
acts on the patch space structures 
$$h|~:~Q''\to Q'$$ for $Q''=Q''_+\cup Q''_-$ 
by acting on the gluing homotopy equivalence 
$$h_{Q''}~:~Q''_0\to Q''_0~.$$ 
Given an element $\sigma\in L_{m-q}(Q'')$ denote its Wall
realization by $\sigma:K\to Q_0''\times [0,1]$, where $\partial
K=Q_0''\amalg \widehat{Q}_0$.  The new patch space structure is given by
$\widehat{Q''}=\widehat{Q''_+}\cup Q''_-$ where
$\widehat{Q''_+}=Q''_+\cup_{id_{Q''_0}} K$.  The new gluing map is now
$h_{\widehat{Q''}}=\sigma|\widehat{Q''}$.  A Poincar\'e bordism from $Q''$ to
$\widehat{Q''}$ is given by $\widehat{Q''_+}\times [0,1] \cup_{\sigma}
Q''_-\times[0,1]$.  where $\sigma$ is viewed as a map $\sigma:K\times
0\to Q_0''\times [0,1]$.
The group $L_m(X'\setminus Q'')$ acts similarly on the structure
$$h_X|~:~X'_0\setminus Q''_0\to X'_0\setminus Q''_0~.$$
The fact that we really had an element of $LS_{m-q}(\Phi)$ insures that
the Wall realizations will glue together to a composite we denote
$K_\sigma$ and that the patch spaces will glue together to give a new
patch space structure, $X''$.
\medskip

\noindent{\bf Claim:} The action of 
$$\sigma\in LS_{m-q}(\Psi)$$ 
results in a new patch space structure whose transversality obstruction
differs by $-\sigma$ from the previous one, i.e.  with the correct
choice of $\sigma$ the transversality obstruction vanishes.
\medskip

\noindent{\bf Proof of claim:} If we consider the addition of $K_\sigma
\times [0,1]$ as in the diagram, it becomes clear that we have changed
$W_0$ by $-K_\sigma$.

\begin{figure}
\epsfysize=1in
\center{\leavevmode \epsfbox{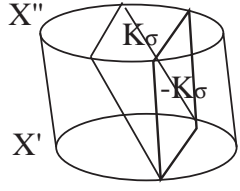}}
\caption{The action of $\sigma$ on $X'$ by $K_\sigma$ produces a new
space $X''$ and adds the above to $W'$.  In particular it adds
$-K_\sigma$ to $W_0$.} \end{figure}

The map $G$ is now homotopic to a Poincar\'e transverse map, and
$G|_X$ is the desired map $g''$ Poincar\'e transverse to $Q$.  The
original map $f$ is homotopic to $g''\circ h$ where $h:M\to X$ is the
homotopy inverse of $H|:X\to M$.\qed
\smallskip

\noindent{\bf Proof of \ref{lem1} (ii):} Suppose given
$$\overline{g}~=~g''\circ h~:~M\to P~,$$ 
for the homotopy equivalence $h:M\to X$ and the Poincar\'e
transverse map $g'':X\to P$.  Let $Q'$ denote $(g'')^{-1}(Q)$, so that by
assumption $Q'$ is a Poincar\'e space.  If $h$ is $s$-cobordant to a
homotopy equivalence $\widehat{h}:\widehat{M}\to X$ which is split over $Q'$ then
$g''\circ \widehat{h}$ is homology manifold transverse to $Q$.  There is an
obstruction $\sigma$ to splitting $h$ which lives in $LS_{m-q}(\Phi)$
where $\Phi$ is defined by the pair $Q'\subset X$.  Let $\sigma$ denote
this obstruction.
\medskip

Because $Q'$ was constructed by the method of \ref{lem1} (i), it has a
two patch structure
$$Q'~=~Q'_+\cup Q'_-\text{ and }Q'_0~=~Q'_+\cap Q'_-$$ 
which agrees with the two patch space structure of $X$, i.e.  $Q'_+\subset
X_+$, $Q'_-\subset X_-$ and $Q'_0\subset X_0$.  If $\Psi$ is defined by
the pair $Q'_0\subset X_0$ there is a natural isomorphism
$$LS_{m-q}(\Psi)~\cong~LS_{m-q}(\Phi)~.$$
$LS_{m-q}(\Psi)$ acts on the possible patch space structures $h:M\to X$
on $M$ as follows.  Recall that by construction, the gluing map
$$h_X~:~X_0\to X_0$$ 
splits along $Q'_0\subset X_0$.  We can change the patch space
structure on $X$, if we make sure that it has a map to the original $X$
which is Poincar\'e transverse to $Q'$.  In particular the new map
$$h_Y~:~Y_0\to Y_0$$ 
must split along the inverse image of $Q'$.  Take $h_X$ and act on it by 
$LS_{m-q}(\Phi)$, by acting on
$$h_{Q'}|~:~Q'_0\to Q'_0$$ 
by $L_{m-q}(Q'_0)$ and by acting on
$$h_X|~:~X_0\setminus Q'_0\to X_0\setminus Q'_0$$ 
by $L_m(X_0\setminus Q'_0)$.  The result fits together along 
$\partial \nu(Q')$ by definition of $LS_{m-q}(\Psi)$.  
By construction, there is a map $k:Y\to X$ which is a homotopy
equivalence, and which is Poincar\'e transverse to $Q'$.
\medskip

We claim that, if $\overline{k}:X\to Y$ is the homotopy inverse of $k$, then
the new splitting obstruction of $$\overline{k}\circ h:M\to X\to Y$$
vanishes.  Clearly the splitting obstruction of $k\circ \overline{k}\circ h$
is the same as that of $h$, but by construction this splitting
obstruction differs from that of $\overline{k}\circ h$ by the $\sigma$
splitting obstruction of $k$.
\medskip

Thus our new map $\overline{k}\circ h$ is $s$-cobordant to a homology
manifold split map.  Composing this homology manifold split map with
the Poincar\'e transverse map $g''\circ k$ we see that our original
$f$ is in fact $s$-cobordant to a homology manifold transverse map.
\qed

\section{Dual transversality  for homology manifolds}

We now extend the results of Chapter 3 on transversality to a codimension $q$
bundle subspace $Q\subset P$ for a map $f:M^m \to P$ from an $m$-dimensional
homology manifold with $m -q \geq 7$, and obtain dual transversality
for a map $f:M^m \to \vert K \vert$ to the polyhedron of a $k$-dimensional
simplicial complex $K$ with $m -k \geq 7$. In Chapter 6 we shall
formulate an obstruction for a map $f:M^m \to \vert K\vert$ to be bordant
to a dual transverse map -- the obstruction is 0 for $m-k \geq 7$.

\medskip

Let $X$ be a simple $m$-dimensional Poincar\'e duality
space with a map $g:X \to \vert K \vert$ to the polyhedron
of a $k$-dimensional simplicial complex $K$.
\medskip

The {\it dual transverse homology manifold structure set}
$\ST^{H,K-\transverse}(X)$ is the set of equivalence classes of pairs
$(M,h)$ with $M$ an $m$-dimensional homology manifold and $h:M \to X$ a
simple homotopy equivalence such that $gh:M \to \vert K \vert$ is dual
transverse, with $(M_1,h_1) \allowbreak \simeq
(M_2,h_2)$ if there exists an $s$-cobordism $(W;M_+,M_2)$ with a simple
homotopy equivalence of the type
$$(f;h_1,h_2)~:~(W;M_1,M_2) \to X \times ([0,1];\{0\},\{1\})$$
such that the composite
$$W \stackrel{f}{\to} X \times [0,1] \stackrel{proj.}{\to} X
\stackrel{g}{\to} \vert K \vert$$
is dual transverse.

\medskip

Let $\TT^{H,K-\transverse}(X)$ be the set of bordism classes of normal maps
$(f,b):M \to X$ from $m$-dimensional homology manifolds such that
$gf:X \to \vert K \vert$ is dual transverse.

\begin{theorem} Let $m-k \geq 7$.\\
{\rm (i)} The dual transverse homology manifold structure set fits into
the surgery exact sequence
$$\dots \to L_{m+1}(\Z[\pi_1(X)]) \stackrel{\partial}{\to}
\ST^{K,Q-\transverse}(X) \stackrel{\eta}{\to}
\TT^{K,Q-\transverse}(X)\stackrel{\theta}{\to} L_m(\Z[\pi_1(X)])~.$$
{\rm (ii)} The forgetful function $\TT^{K,Q-\transverse}(X) \to \TT^K(X)$
is a bijection.\\
{\rm (iii)} The forgetful function $\ST^{K,Q-\transverse}(X) \to \ST^K(X)$
is a bijection. In particular, $\ST^{K,Q-\transverse}(X)$ is non-empty if and only if
$\ST^K(X)$ is non-empty.\\
{\rm (iv)} Every map $f:M \to \vert K \vert$ from an $m$-dimensional homology
manifold is dual $s$-transverse.
\end{theorem}
\begin{proof} Exactly as above.
\end{proof}

\section{The dual transversality obstruction}

For any simplicial complex $K$ let $\Omega^H_m(K)$ (resp.
$\Omega^{H,\transverse}_m(K)$) be the bordism group of maps $f:M \to
\vert K \vert$ from $m$-dimensional homology manifolds (resp.  dual
transverse maps $f$).
We now formulate the bordism obstruction to dual transversality for homology
manifolds, as the failure of the forgetful map
$\Omega^{H,\transverse}_m(K) \to \Omega_m^H(K)$ to be an isomorphism.
\medskip

We refer to Chapters 11,12 of Ranicki \cite{R3} for an exposition of
$\Delta$-sets, generalized homology theories, bordism spectra and
assembly maps.  The topological manifold bordism groups and the
$\Omega^{TOP}_{\bullet} (\{ *\})$-coefficient generalized homology
groups are the homotopy groups of spectra of Kan $\Delta$-sets
$\Omega^{TOP}_{\bullet} (K)$, 
$\HH_{\bullet}(K;\Omega^{TOP}_{\bullet}(\{ *\}))$
$$\begin{array}{l}
\Omega^{TOP}_*(K)~=~\pi_*(\Omega^{TOP}_{\bullet} (K))~,\\[1ex]
H_*(K;\Omega^{TOP}_{\bullet} (\{ *\}))~=~
\pi_*(\HH_{\bullet}(K;\Omega^{TOP}_{\bullet} (\{ *\})))~.
\end{array}$$
Moreover, there is defined a topological assembly map 
$$A^{TOP}~:~\HH_{\bullet}(K;\Omega^{TOP}_{\bullet} (\{ *\})) 
\to \Omega^{TOP}_{\bullet} (K)$$
which is a homotopy equivalence by topological manifold transversality, 
inducing the Pontrjagin-Thom isomorphisms
$$A^{TOP}~: ~H_*(K;\Omega^{TOP}_{\bullet} (\{ *\}))~\cong~\Omega^{TOP}_*(K)$$
The combinatorial construction of $\Omega^{TOP}_{\bullet} (K)$ and $A^{TOP}$
will now be extended to the homology manifold bordism spectra $\Omega^H_{\bullet} (K)$,
$\Omega^{H,\transverse}_{\bullet} (K)$ with an assembly map
$$A^H~: ~\Omega^{H,\transverse}_{\bullet} (K)~=~
\HH_{\bullet}(K;\Omega^H_{\bullet} (\{ * \})) \to \Omega^H_{\bullet} (K)~.$$

The homology assembly map $A^H$ is only a homotopy equivalence
to the extent to which homology manifolds have transversality.
\medskip

The {\it homology manifold bordism spectrum}
of a simplicial complex $K$
$$\Omega^H_{\bullet} (K)=\{\Omega^H_{\bullet} (K)_m\,\vert\,m\geq 0\}
\vspace*{1mm}$$
is the $\Omega$-spectrum with $\Omega^H_{\bullet} (K)_m$ the Kan
$\Delta$-set defined by
$$\begin{array}{rl}
\Omega^H_{\bullet} (K)^{(n)}_m~=~\{\,\!\!\!\!\!
&\hbox{$(m+n)$-dimensional homology manifold $n$-ads}
\vspace*{1mm}\\
&(M;\partial_0M,\partial_1M,\dots,\partial_nM) \hbox{ such that }
\vspace*{1mm}\\
&\partial_0M \cap \partial_1M \cap \dots \cap \partial_nM=\emptyset,
\,\hbox{\rm with a map } f:M \to \vert K \vert~\}
\end{array}$$
with base simplex the empty manifold $n$-ad $\emptyset$.
The homotopy groups
$$\pi_m(\Omega^H_{\bullet} (K))~=~\pi_{m-k}(\Omega^H_{\bullet} (K)_k)~=~
\Omega^H_m(K)~~(m\geq k \geq 0)$$
are the bordism groups of maps $f:M \to \vert K \vert$ from closed
$m$-dimensional homology manifolds.  Similarly for the {\it dual
transverse bordism spectrum} $\Omega^{H,\transverse}_{\bullet} (K)$,
with the additional requirement that $f:M \to \vert K \vert$ be dual
transverse.
\medskip

\begin{prop} \label{four}
The dual transverse bordism spectrum $\Omega^{H,\transverse}_*(K)$
coincides with the generalized $\Omega^H_{\bullet} (\{ * \})$-homology
spectrum of $K$
$$\Omega^{H,\transverse}_{\bullet}(K)~=~\HH(K;\Omega^H_{\bullet} (\{ * \}))~,$$
so that on the level of homotopy groups
$$\Omega^{H,\transverse}_m(K)~=~H_m(K;\Omega^H_{\bullet} (\{ * \}))~~(m \geq 0)~.$$
\end{prop}
\begin{proof}

Define the generalized homology spectrum
$\HH_{\bullet}(K;\Omega^H_{\bullet} (\{ * \}))$ as in
\cite{R3} (12.3), with
an $m$-dimensional $\Omega^H_{\bullet} (\{*\})$-coefficient cycle
(\cite{R3},12.17)
$$x~=~\{M(\sigma )^{m- \vert \sigma \vert}\, \vert \, \sigma \in K \}$$
essentially the same as a dual transverse map
$$f(x)~:~ M^m~=~\bigcup_{\sigma \in K}M(\sigma ) \to \vert K \vert$$
from a closed $m$-dimensional homology manifold, with inverse images
the $(m-\vert \sigma \vert)$-dimensional homology manifolds with 
$$f(x)^{-1}(D(\sigma ,K),\partial D(\sigma ,K))~=~
(M(\sigma),\partial M(\sigma ))~~(\sigma \in K)~.$$
The homotopy group
$$\pi_m(\HH_{\bullet}(K;\Omega^H_{\bullet} (\{ * \})))~=~
H_m(K;\Omega^H_{\bullet} (\{ * \}))$$
is the cobordism group of such cycles, and is the bordism group of dual
transverse maps $f: M \to \vert K \vert$ from $m$-dimensional homology
manifolds.
\end{proof}

The {\it homology assembly map} of spectra
$$A^H~: ~\Omega^{H,\transverse}_{\bullet}(K)~=~
\HH_{\bullet}(K;\Omega^H_{\bullet} (\{ * \})) \to \Omega^H_{\bullet} (K)$$
is defined as in \cite{R3} (12.18), inducing on the level of homotopy groups
the assembly maps of bordism groups
$$A^H~:~\Omega^{H,\transverse}_*(K)~=~H_*(K;\Omega^H_{\bullet} (\{ * \})) \to
\Omega^H_*(K)$$
which forget dual transversality.

\medskip

\begin{definition} {\normalfont
Given an $m$-dimensional homology manifold
$$M~=~\bigcup_{\alpha}M_{\alpha}$$
with $M_{\alpha}$ the components of $M$, set
$$\I(M)~=~\sum\limits_{\alpha}[M_{\alpha}] i(M_{\alpha})\in H_m(M)[L_0(\Z)]~,$$
with $[M_{\alpha}]\in H_m(M)$ the image of the fundamental class
$[M_{\alpha}]\in H_m(M_{\alpha})$.}
\end{definition}

The augmentation map
$$H_m(M)[L_0(\Z)] \to H_m(M;L_0(\Z))~=~H_m(M)\otimes_{\Z}L_0(\Z)~;~x[y] \to x\otimes y$$
sends $\I(M)$ to the resolution obstruction of $M$
$$i(M)~=~\sum\limits_{\alpha}i(M_{\alpha}) \in H_m(M;L_0(\Z))~=~\sum\limits_{\alpha}
L_0(\Z)~.$$

We shall now verify that $\I(M,f)$ and $\widetilde\I(M,f)$ are
bordism invariants\,:

\begin{prop} \label{bord}
Let $f:M \to \vert K \vert$ be a map from an $m$-dimensional homology manifold
$M$ to a polyhedron, and let $M=\bigcup\limits_{\alpha}M_{\alpha}$ be the decomposition
of $M$ into components with $f_\alpha=f\vert:M_{\alpha} \to \vert
K\vert$.\\
{\rm (i)} The element
$$\I(M,f)~=~\sum\limits_{\alpha}(f_{\alpha})_*[M_{\alpha}][i(M_{\alpha})]
\in H_m(K)[L_0(\Z)]$$
is a homology manifold bordism invariant.\\
{\rm (ii)} If $f$ is a dual transverse map then
$$\I(M,f)~=~0 \in H_m(K)[L_0(\Z)]~.$$
{\rm (iii)} If $f_*:H_m(M) \to H_m(K)$ is an isomorphism
and $f$ is homology manifold bordant to a dual transverse map then
$$i(M)~=~0 \in H_m(M;L_0(\Z))~.$$
\end{prop}
\begin{proof} {\rm (i)} Given a bordism
$$(g;f,f')~:~(W;M,M') \to \vert K \vert~.$$
Denote the connected components by
$$(g_\alpha;f_\alpha,f'_\alpha)~:~(W_\alpha;M_\alpha,M'_\alpha) \to \vert K \vert.$$
We have
$$(f_\alpha)_*[M_\alpha]~=~(f'_\alpha)_*[M'_\alpha]\in H_m(K)$$ as usual.
\medskip

Furthermore $$i(M_\alpha)=i(W_\alpha)=i(M'_\alpha)$$  follows from the fact that
$i(X)=i(\partial X)$ for any connected homology manifold with non-empty boundary
(Quinn \cite{Q2}, 1.1).  Thus we have
$$\begin{array}{ll}
\I(M,f)\!\!&=~
\sum\limits_{\alpha}(f_{\alpha})_*[M_{\alpha}][i(M_{\alpha})]\\
&=~\sum\limits_{\alpha'}(f'_{\alpha'})_*[i(M_{\alpha'})]~
=~ \I(M',f')\in H_m(K)[L_0(\Z)]~.
\end{array}
$$
{\rm (ii)} The augmentation map
$$H_m(K)[L_0(\Z)] \to H_m(K;L_0(\Z))~=~H_m(K)\otimes_{\Z}L_0(\Z)~;~x[y] \to x\otimes y$$
sends $\I(M,f)$ to $f_*i(M) \in H_m(K;L_0(\Z))$, for any map
$f:M \to \vert K \vert $. If $M$ is connected then $i(M) \in L_0(\Z)$ and
$$\I(M,f) ~=~f_*[M][i(M)] \in H_m(K)[i(M)] \subset H_m(K)[L_0(\Z)]~,$$
so that $f_*i(M)=0$ implies $\I(M,f) =0$.
If $f:M \to \vert K \vert $ is dual transverse
$$f_*i(M)~=~\sum\limits_{\sigma \in K ^{(m)}}i(f^{-1}D(\sigma,K)) \sigma
\in H_m(K;L_0(\Z))~.$$
Each $f^{-1}D(\sigma,K)$ is a 0-dimensional homology manifold,
which is a disjoint union of points, with
resolution obstruction $i(f^{-1}D(\sigma,K))=0$, so that $f_*i(M)=0$.
Thus if $M$ is connected and $f:M \to \vert K \vert $ is dual transverse
then $f_*i(M)=0$, and hence $\I(M,f)=0$. Apply this to each
component of $M$.\\
{\rm (iii)} Combine {\rm (i)} and {\rm (ii)}.
\end{proof}

\begin{definition} {\normalfont
{\rm (i)} For any space $X$ define a morphism
$$\I~:~\Omega^H_m(X) \to H_m(X)[L_0(\Z)]~;~ (M,f) \mapsto
\I(M,f)=f_*\I(M)$$
such that
$$i~:~\Omega^H_m(X)
\stackrel\I{\to} H_m(X)[L_0(\Z)]  \to H_m(X;L_0(\Z))~;~
(M,f) \mapsto f_*i(M)$$
and
$$\Omega^H_m(X) \stackrel\I{\to} H_m(X)[L_0(\Z)]
\stackrel{proj.}{\longrightarrow} H_m(X)~;~ (M,f) \mapsto f_*[M]~.$$
{\rm (ii)} Let
$$\widetilde{L}_0(\Z)~=~L_0(\Z)\backslash \{0\}~,$$
and let
$$\widetilde\I~:~\Omega^H_m(X)\stackrel\I{\to}
H_m(X)[L_0(\Z)] \to H_m(X)[\widetilde{L}_0(\Z)]~.$$}
\end{definition}

For any $\Omega$-spectrum $X_{\bullet}=\{X_n\,\vert\,n \geq 0\}$
of $\Delta$-sets there is defined an $\Omega$-spectrum of $\Delta$-sets
$$X_{\bullet}[L_0(\Z)]~=~\{X_n[L_0(\Z)]\,\vert\,n \geq 0\}$$
with
$$\pi_m(X[L_0(\Z)])~=~\pi_m(X)[L_0(\Z)]~~(m \geq 0)~.$$

\begin{prop} \label{five}
Let $K$ be a $k$-dimensional simplicial complex.\\
{\rm (i)} The map of $\Omega$-spectra
$$\Omega^H_{\bullet}(K) \to \Omega^{TOP}_{\bullet}(K)[L_0(\Z)]~;~M \mapsto M_{TOP}[i(M)]$$
induces maps
$$\pi_m(\Omega_{\bullet}^H(K))=\Omega^H_m(K) \to
\pi_m(\Omega^{TOP}_{\bullet}(K)[L_0(\Z)])=\Omega_m^{TOP}(K)[L_0(\Z)]$$
which are isomorphisms for $m \geq 6$.\\
{\rm (ii)} The composite
$$\Omega^{H,\transverse}_m(K) \stackrel{A^H}{\to}
\Omega^H_m(K) \stackrel{\widetilde\I}{\to}  H_m(K)[\widetilde{L}_0(\Z)]$$
is 0.\\
{\rm (iii)} The homology assembly map
$$A^H~:~\Omega^{H,\transverse}_m(K)~=~H_m(K;\Omega^H_{\bullet} (\{ * \}))
\to \Omega^H_m(K)$$
is an isomorphism for $m-k\geq 6$.
\end{prop}
\begin{proof} {\rm (i)} See Johnston \cite{Jo} (cf. Corollary \ref{jbord}).\\
{\rm (ii)} By Proposition \ref{bord}
$$\text{im}(\I A^H:\Omega^{H,\transverse}_m(K) \to H_m(K)[L_0(\Z)])
\subseteq H_m(K)[0] \subseteq H_m(K)[L_0(\Z)]~,$$
so that $\widetilde\I A^H=0$.\\
{\rm (iii)} By {\rm (i)} we have that the maps of Kan $\Delta$-sets
$$\Omega^H_{\bullet}(\{*\})_n \to \Omega^{TOP}_{\bullet}(\{*\})_n[L_0(\Z)]~;~
M \mapsto M_{TOP}[i(M)]$$
are homotopy equivalences for $n \geq 6$.
It follows that for $m \geq k+6$
$$\HH(K;\Omega^H_{\bullet}(\{*\}))_m \simeq
\HH(K;\Omega^{TOP}_{\bullet}(\{*\}))_m[L_0(\Z)]$$
and hence that
$$\begin{array}{rl}
H_m(K;\Omega^H_{\bullet} (\{ * \}))\!\!
&=~\pi_0(\HH(K;\Omega^H_{\bullet}(\{*\}))_m)\vspace*{1mm}\\
&=~\pi_0(\HH(K;\Omega^{TOP}_{\bullet}(\{*\}))_m[L_0(\Z)])\vspace*{1mm}\\
&=~H_m(K;\Omega^{TOP}_{\bullet}(\{*\}))[L_0(\Z)]\vspace*{1mm}\\
&=~\Omega^{TOP}_m(K)[L_0(\Z)]~=~\Omega^H_m(K)~.
\end{array}$$
\end{proof}

\begin{theorem}\label{two}
{\rm (i)}
The composite
$$\Omega^{H,\transverse}_m(K) \to \Omega^H_m(K)
\stackrel{\widetilde\I}{\to} H_m(K)[\widetilde{L}_0(\Z)]$$
is 0. Thus if $f:M \to \vert K\vert$ is homology manifold bordant
to a dual transverse map then
$$\I(M,f) \in H_m(K)[0] \subset H_m(K)[L_0(\Z)]$$
and either $f_*[M]=0 \in H_m(K)$ or $i(M)=0 \in H_m(M;L_0(\Z))$,
but in any case
$$f_*i(M)=0 \in H_m(K;L_0(\Z))~.$$
{\rm (ii)} There exists a spectrum $\LLL_{\bullet}$ whose homotopy groups 
fit into an exact sequence
$$\dots \to\Omega^H_m(\{*\}) \to \Omega^{TOP}_m(\{*\})[L_0(\Z)] \to
\pi_m(\LLL_{\bullet})\to \Omega^H_{m-1}(\{*\}) \to \dots~,$$
with
$$\pi_m(\LLL_{\bullet})~=~
\begin{cases}
\Z[\widetilde{L}_0(\Z)]&\text{if $m=0$}\\
0&\text{if $m\geq 1$ and $m \neq 4,5$}
\end{cases}$$
and
$$\begin{array}{ll}
&\Omega^H_m(K)\to \Omega^{TOP}_m(K)[L_0(\Z)] \stackrel{\widetilde\I}{\to}
H_m(K;\LLL_{\bullet})\vspace*{2mm}\\
&\hphantom{\Omega^H_m(K)\to \Omega^{TOP}_m(K)[L_0(\Z)]}
\to H_m(K;\pi_0(\LLL_{\bullet}))~=~H_m(K)[\widetilde{L}_0(\Z)]~;\vspace*{2mm}\\
&\hphantom{\Omega^H_m(K)\to \Omega^{TOP}_m(K)[L_0(\Z)]}
(M,f:M \to \vert K \vert) \mapsto f_*[M]\widetilde{[i(M)]}~.\\
\end{array}$$
{\rm (iii)} For $m \geq 6$ there is defined an exact sequence
$$\dots \to\Omega^{H,\transverse}_m(K) \stackrel{A^H}{\to} \Omega^H_m(K)
\stackrel{\widetilde\I}{\to}
H_m(K;\LLL_{\bullet})\to \Omega^{H,\transverse}_{m-1}(K) \to \dots~.$$
\end{theorem}

The proof of \ref{two} {\rm (i)} is given by Proposition \ref{bord}.
The remainder of this Chapter is devoted to the proofs of \ref{two} {\rm (ii)}
and {\rm (iii)}. The spectrum $\LLL_{\bullet}$ in {\rm (ii)} is given by\,:

\begin{definition}{\normalfont Let
$$\LLL_{\bullet}~=~
\text{cofibre}(\Omega^H_{\bullet}(\{*\}) \to
\Omega^{TOP}_{\bullet}(\{*\})[L_0(\Z)])~,$$
an $\Omega$-spectrum whose homotopy groups fit into an exact sequence
$$\dots \to\Omega^H_m(\{*\}) \to \Omega^{TOP}_m(\{*\})[L_0(\Z)] \to
\pi_m(\LLL_{\bullet})\to \Omega^H_{m-1}(\{*\}) \to \dots~.$$}
\end{definition}
\medskip

It is now immediate from the identities
$$\begin{array}{l}
\Omega^H_m(\{*\})~=~\begin{cases}
\Z&\text{if $m=0$} \\
0&\text{if $m=1,2$}\\
\Omega^{TOP}_m(\{*\})[L_0(\Z)]&\text{if $m\geq 6$~,}
\end{cases}\vspace*{2mm}\\
\Omega^{TOP}_m(\{*\})~=~\begin{cases}
\Z&\text{if $m=0$} \\
0&\text{if $m=1,2,3$}
\end{cases}
\end{array}$$
that
$$\pi_m(\LLL_{\bullet})~=~\begin{cases}
\Z[\widetilde{L}_0(\Z)]&\text{if $m=0$}\\
0&\text{if $m\geq 1$ and $m \neq 4,5$}
\end{cases}$$
as claimed in the statement of Theorem \ref{two} {\rm (ii)}.
\medskip

The exact sequence in the statement of Theorem \ref{two} {\rm (iii)} is given by\,:

\begin{prop}\label{six} For $m \geq 6$ there is defined an exact sequence
$$\dots \to\Omega^{H,\transverse}_m(K)
\stackrel{A^H}{\to} \Omega^H_m(K)
\stackrel{\widetilde\I}{\to}
H_m(K;\LLL_{\bullet})\to \Omega^{H,\transverse}_{m-1}(K) \to \dots~.$$
\end{prop}
\begin{proof} This is just the exact sequence
$$\begin{array}{l}
\dots \to H_m(K;\Omega_{\bullet}^H(\{*\})) \to
H_m(K;\Omega_{\bullet}^{TOP}(\{*\})[L_0(\Z)])\vspace*{2mm}\\
\hphantom{\dots \to H_m(K;\Omega_{\bullet}^H(\{*\}))}
\to H_m(K;\LLL_{\bullet})\to H_{m-1}(K;\Omega_{\bullet}^H(\{*\})) \to \dots
\end{array}$$
induced by the cofibration sequence of spectra
$$\Omega_{\bullet}^H(\{*\}) \to  \Omega_{\bullet}^{TOP}(\{*\})[L_0(\Z)]
\to \LLL_{\bullet}~,$$
using \ref{jbord}, \ref{four} to identify
$$\begin{array}{l}
\Omega^{H,\transverse}_m(K)~=~H_m(K;\Omega_{\bullet}^H(\{*\}))~,\vspace*{2mm}\\
\Omega^H_m(K)~=~H_m(K;\Omega_{\bullet}^{TOP}(\{*\}))[L_0(\Z)]~.
\end{array}$$
\end{proof}

This completes the proof of Theorem \ref{two}.

\medskip

In the special case $K=S^m$ we have\,:

\begin{cor}{
{\rm (i)} A map $f:M^m \to S^m$ from an $m$-dimensional homology manifold $M$
determines the element
$$\widetilde\I(M,f)~=~{\rm degree}(f)[i(M)] \in \Z[\widetilde{L}_0(\Z)]~,$$
which vanishes if (and for $m\geq 6$ only if) $f$ is bordant to a dual transverse map.\\
{\rm (ii)} For connected $M$ $\widetilde\I(M,f)=0$ if and only if
either ${\rm degree}(f)=0 \in \Z$ or $i(M)=0 \in H_m(M;L_0(\Z))=L_0(\Z)$.}
\end{cor}
\begin{proof} This follows from the sequence of \ref{two} (i)
$$ H_m(S^m; \Omega^H_{\bullet}(\{*\}))
\stackrel{A^H}{\longrightarrow}
\Omega^H_m(S^m)\stackrel{\widetilde\I}{\longrightarrow}\Z[\widetilde{L}_0(\Z)]$$ whose
composite is 0, and which is exact for $m \geq 6$.
By a simple calculation
$$H_m(S^m; \Omega^H_{\bullet}(\{*\})) ~=~
\Omega^H_0(\{*\}) \oplus \Omega^H_m(\{*\})~.$$
Since any
zero-dimensional homology manifold is a disjoint union of points
$$\Omega^H_0(\{*\})~=~ \Omega^{TOP}_0(\{*\})~=~ \Z~.$$
By \ref{jbord}
$$\Omega^H_m(\{*\}) ~=~ \Omega_m^{TOP}(\{*\})[L_0(\Z)]~~(m\geq 6)~,$$
so that
$$H_m(S^m; \Omega^H_{\bullet}(\{*\})) ~=~
\Omega^{TOP}_0(\{*\}) \oplus \Omega^{TOP}_m(\{*\})[L_0(\Z)]~.$$
Also
$$\Omega^H_m(S^m)~=~ \Omega^{TOP}_m(S^m)[L_0(\Z)]~~(m\geq 6)~,$$
and by topological transversality
$$\Omega^{TOP}_m(S^m)~=~ H_m(S^m; \Omega^{TOP}_{\bullet}(\{*\}))
~=~ \Omega^{TOP}_0(\{*\}) \oplus
\Omega^{TOP}_m(\{*\})~,$$
so that
$$\Omega^H_m(S^m)~=~ \Omega^{TOP}_0(\{*\})[L_0(\Z)] \oplus
\Omega^{TOP}_m(\{*\})[L_0(\Z)]~~(m\geq 6)~.$$
The components of the assembly map
$$\begin{array}{l}
A^H~=~A^H_0\oplus A^H_m~:\\[1ex]
\Omega^H_0(\{*\}) \oplus \Omega^H_m(\{*\}) \to \Omega^{TOP}_0(\{*\})[L_0(\Z)]
\oplus \Omega^{TOP}_m(\{*\})[L_0(\Z)]
\end{array}$$
are given by the inclusion
$$A^H_0~:~\Omega^H_0(\{*\})~=~\Z\to
\Omega^{TOP}_0(\{*\})[L_0(\Z)]~=~\Z[L_0(\Z)]~;~t \mapsto t[0]$$
and the isomorphism
$$A^H_m~:~\Omega^H_m(\{*\})~\cong~ \Omega^{TOP}_m(\{*\})[L_0(\Z)]~.$$
The cokernel of $A^H$ is thus given by the cokernel of $A^H_0$ as
$\Z[\widetilde{L}_0(\Z)]$.  In particular, if $m\geq 6$ and $\Sigma^m$
is one of the $m$-dimensional homology spheres with 
$$i(\Sigma^m) \neq 0 \in H_m(\Sigma^m;L_0(\Z))~=~L_0(\Z)$$ 
constructed in \cite{BFMW} there exists a homotopy equivalence
$f:\Sigma^m \to S^m$ with
$$\widetilde\I(\Sigma^m,f)~=~i(\Sigma^m) \neq 0 \in \Z[\widetilde{L}_0(\Z)]~,$$
so that $f$ is not bordant to a dual transverse map.
\end{proof}

Every compact $ANR$ $X$ is simple homotopy equivalent to a polyhedron
$\vert K\vert$, by the result of West \cite{Wes}.  However, if $X$ is an
$m$-dimensional homology manifold then $\vert K\vert$ may not be a
homology manifold.
\medskip

\noindent{\bf Theorem} 
{\rm (Levitt and Ranicki \cite{LR}, Ranicki \cite{R3} (19.6))}\newline
{\it Let $X$ be a compact $ANR$ which is a simple $m$-dimensional Poincar\'e
space, and let $f:X \to \vert K\vert$ be a simple homotopy equivalence
to a polyhedron.  The map $f$ is simple Poincar\'e bordant to a dual
Poincar\'e transverse map $f':X' \to \vert K \vert$ if {\rm (}and for
$m\geq 6$ only if\/{\rm )} $X$ is simple homotopy equivalent to a
topological manifold.} 
\medskip

Theorem \ref{two} will now be used to obtain an analogous characterization
of resolvable homology manifolds.

\begin{prop} \label{three}
Let $M$ be an $m$-dimensional homology manifold, and let $f:M \to \vert
K\vert$ be a simple homotopy equivalence to a polyhedron.  The map $f$
is $s$-cobordant to a dual transverse homotopy equivalence $f':M' \to
\vert K \vert$ if {\rm (}and for $m\geq 6$ only if\/{\rm )} $M$ is resolvable.
\end{prop}
\begin{proof}
If $M$ is resolvable the mapping cylinder of a resolution $M_{TOP} \to M$ is an
$s$-cobordism $(W;M,M_{TOP})$ with a bordism
$$(g;f,f_{TOP})~:~(W;M,M_{TOP}) \to \vert K \vert$$
to a homotopy equivalence
$f_{TOP}:M_{TOP} \to \vert K\vert$ which is topologically dual transverse.\\
\indent Conversely, suppose that $f:M \to \vert K \vert$ is a dual transverse
homotopy equivalence. Without loss of generality, it may be assumed that
$M$ is connected, so that
$$\I(M)~=~[i(M)] \in \Z[i(M)] \subseteq H_m(M)[L_0(\Z)]~=~\Z[L_0(\Z)]~.$$
The composite
$$\Omega^H_m(K) \stackrel{\widetilde\I}{\to}
H_m(K;\LLL_{\bullet}) \to H_m(K)[\widetilde{L}_0(\Z)]$$
sends $(M,f)$ to
$$\widetilde\I(M,f)~=~
\begin{cases}
[i(M)] \in \Z[i(M)] \subseteq \Z[\widetilde{L}_0(\Z)]&\text{if $i(M)\neq 0$}\\
0 \in \Z[\widetilde{L}_0(\Z)]&\text{if $i(M)=0$~.}
\end{cases}
$$
Since $f:M \to \vert K \vert$ is dual transverse,
$i(M)=0 \in L_0(\Z)$ by \ref{two} (ii).
\end{proof}

\section{The multiplicative structure of homology manifold bordism}

The product of an $m$-dimensional Poin\-car\'e space $X$ and
an $n$-dimensional Poincar\'e space $Y$ is an $(m+n)$-dimensional
Poincar\'e space $X \times Y$, with Spivak normal fibration
$$\nu_{X\times Y} ~=~\nu_X\times \nu_Y ~:~X\times Y \to BG~.$$
If $X$ is an $m$-dimensional homology manifold and $Y$ is an $n$-dimensional
homology manifold then $X \times Y$ is an $(m+n)$-dimensional homology
manifold, with the resolution obstructions satisfying the index
product formula of Quinn \cite{Q2}
$$1+8\,i(X\times Y)~=~(1+8\,i(X))(1+8\,i(Y)) \in 1+8\Z~.$$
In general, the canonical $TOP$ reductions of the Spivak normal fibrations
of $X,Y,$\\ $X\times Y$ are such that
$$\widetilde{\nu}_{X\times Y} ~\neq~
\widetilde{\nu}_X\times \widetilde{\nu}_Y ~:~X\times Y \to BTOP~.$$
We shall now analyze the difference
$$\widetilde{\nu}_{X\times Y} -\widetilde{\nu}_X\times \widetilde{\nu}_Y~:~
X \times Y \to G/TOP$$
using the multiplicative properties of the algebraic $L$-spectra of \cite{R3}.
This analysis will be used to obtain the product
structure on $\Omega^{TOP}_*(K)[L_0(\Z)]$ which corresponds to the
product structure on $\Omega^H_*$
$$\Omega_m^H(J) \otimes \Omega^H_n(K) \to \Omega_{m+n}^H(J\times K)~;~
X\otimes Y \mapsto X \times Y$$
under the isomorphism of abelian groups given by \ref{jbord}
$$\Omega^H_*(K) \stackrel{\cong}{\to} \Omega^{TOP}_*(K)[L_0(\Z)]~;~
X \mapsto X_{TOP}[i(X)]~.$$

\medskip

The Spivak normal fibration $\nu_X:X \to BG(k)$ ($k$ large)
of an $m$-dimensional Poincar\'e space $X$ is equipped with a degree 1 map
$\rho:S^{m+k} \to T(\nu_X)$. If $\nu_X$ is $TOP$ reducible
then for any $TOP$ reduction $\widetilde{\nu}_X:X \to BTOP$
the Browder-Novikov transversality construction gives a degree 1 normal map
$$f~=~\rho\vert~:~M~=~\rho^{-1}(X) \to X$$
with $M$ an $m$-dimensional topological manifold.
If $\widetilde{\nu},\widetilde{\nu}':X \to BTOP$ are two $TOP$ reductions
of $\nu_X$ the difference is classified by an element
$$t(\widetilde{\nu},\widetilde{\nu}') \in H_m(X;\LL_{\bullet})~=~[X,G/TOP]$$
and the corresponding degree 1 normal maps $f:M \to X$, $f':M' \to X$
can be chosen to be such that $f'=fg$ for a degree 1 normal map
$g:M' \to M$ classified by an element
$$[g]_{\LL} \in H_m(M;\LL_{\bullet})~=~[M,G/TOP]$$
such that $f_*[g]_{\LL}=t(\widetilde{\nu},\widetilde{\nu}')$.
The surgery obstruction of $g$ is the assembly of $[g]_{\LL}$
$$\sigma_*(g)~=~A([g]_{\LL}) \in L_m(\Z[\pi_1(M)])$$
and the surgery obstructions of $f,f'$ differ by
$$\sigma_*(f') - \sigma_*(f)~=~A(f_*[g]_{\LL})
~=~A(t(\widetilde{\nu},\widetilde{\nu}')) \in L_m(\Z[\pi_1(X)])~.$$
See Chapter 16 of \cite{R3} for the $L$-theory orientation of topology.

\medskip

Let $\LL^{\bullet}$ be the 0-connective symmetric ring $L$-spectrum of $\Z$,
with homotopy groups
$$\pi_m(\LL^{\bullet})~=~L^m(\Z)~=~
\begin{cases}
\Z&\text{if $m \equiv 0$ (mod 4)}\\
\Z_2&\text{if $m \equiv 1$ (mod 4)}\\
0&\text{if $m \equiv 2$ (mod 4)}\\
0&\text{if $m \equiv 3$ (mod 4)~.}\\
\end{cases}$$

\begin{theorem}\label{sixteen} {\rm (Ranicki \cite{R3}, 25.7)}\newline
An $m$-dimensional homology manifold $X$ has a
canonical $\LL^{\bullet}$-orientation
$$[X]_{\LL}\in H_m(X;\LL^{\bullet})$$
with assembly the symmetric signature of $X$
$$A([X]_{\LL}) ~=~\sigma^*(X) \in L^m(\Z[\pi_1(X)])~,$$
and such that there are defined $\LL^{\bullet}$-coefficient Poincar\'e duality isomorphisms
$$[X]_{\LL} \cap - ~:~ H^*(X;\LL^{\bullet}) \cong H_{m-*}(X;\LL^{\bullet})~,$$
as well as with coefficients in any $\LL^{\bullet}$-module spectrum
{\rm (}e.g. $\LL_{\bullet}$, $\overline{\LL}_{\bullet}$\/{\rm )}.
\end{theorem}

\begin{proof} See \cite{R3} (16.16) for the canonical
$\LL^{\bullet}$-orientation of an $m$-dimensional topological manifold.
Let $f:M=X_{TOP}\to X$ be the normal map from a topological manifold
determined (up to normal bordism) by the canonical $TOP$ reduction
$\widetilde{\nu}_X$ of $\nu_X$, with surgery obstruction

$$\sigma_*(f)~=~ \overline{A}([f]_{\LL})~=~\overline{A}(-i(X)) \in L_m(\Z[\pi_1(X)])$$
the assembly of
$$[f]_{\LL}~=~(-i(X),0) \in H_m(X;\overline{\LL}_{\bullet})~=~
H_m(X;L_0(\Z)) \oplus H_m(X;\LL_{\bullet})~.$$
The canonical $\LL^{\bullet}$-theory orientation of $X$ is given by
$$\begin{array}{ll}
[X]_{\LL}\!\!\!
&=~f_*[M]_{\LL} -(1+T)[f]_{\LL}~=~f_*[M]_{\LL} +(8\,i(X),0) \vspace*{2mm}\\
&\in H_m(X;\LL^{\bullet})~=~
H_m(X;L^0(\Z)) \oplus H_m(X;\LL^{\bullet}\langle 1 \rangle)~.
\end{array}$$
\end{proof}

Let $\LL^{\bullet}\langle 1 \rangle$ be the 1-connective cover
of $\LL^{\bullet}$, so that for any space $X$ there is defined an exact sequence
$$\dots \to H_m(X;\LL^{\bullet}\langle 1 \rangle) \to H_m(X;\LL^{\bullet})
\to H_m(X;L^0(\Z)) \to H_{m-1}(X;\LL^{\bullet}\langle 1 \rangle)\to \dots~.$$
The 0th spaces of $\LL^{\bullet}$ and $\LL^{\bullet}\langle 1 \rangle$
are related by
$$\LL^0~=~L^0(\Z) \times \LL^0\langle 1 \rangle~,$$
so that for an $m$-dimensional homology manifold $X$
$$\begin{array}{ll}
H_m(X;\LL^{\bullet})
&\cong~ H^0(X;\LL^{\bullet})\vspace*{2mm}\\
&\cong~ [X,\LL^0]~\cong~ [X,L^0(\Z) \times \LL^0\langle 1 \rangle]\vspace*{2mm}\\
&\cong~ H^0(X;L^0(\Z)) \oplus H^0(X;\LL^{\bullet}\langle 1 \rangle)\vspace*{2mm} \\
&\cong~ H_m(X;L^0(\Z)) \oplus H_m(X;\LL^{\bullet}\langle 1 \rangle)~.
\end{array}$$
For a connected $X$
$$[X]_{\LL}~=~(1+8\,i(X),\widetilde{[X]}_{\LL})\in H_m(X;\LL^{\bullet})~=~
H_m(X;L^0(\Z)) \oplus H_m(X;\LL^{\bullet}\langle 1 \rangle)~,$$
writing the ordinary fundamental class as $[X]=1\in H_m(X;L^0(\Z))=\Z$.

\begin{cor}
Given an $m$-dimensional homology manifold $X$ let
$f:M=X_{TOP}\allowbreak \to X$ be a normal map from a topological
manifold in the canonical class.  The canonical $\LL^{\bullet}$-orientation 
of $M$ is such that
$$[M]_{\LL}~=~(1,\widetilde{[M]}_{\LL})\in H_m(M;\LL^{\bullet})~=~
H_m(M;L^0(\Z)) \oplus H_m(M;\LL^{\bullet}\langle 1 \rangle)$$
with
$$\begin{array}{l}
f_*[M]_{\LL}~=~-8\,i(X)+[X]_{\LL}\in  H_m(X;\LL^{\bullet})~,\vspace*{2mm}\\
f_*\widetilde{[M]}_{\LL}~=~
\widetilde{[X]}_{\LL} \in  H_m(X;\LL^{\bullet}\langle 1 \rangle)~.
\end{array}
$$
\end{cor}

\begin{definition} {\normalfont
Let $X$ be an $m$-dimensional homology manifold, and let
$\widetilde{\nu}_X:X \to BTOP$ be the canonical Ferry-Pedersen
\cite{FP} $TOP$ reduction of the Spivak normal fibration $\nu_X:X \to BG$, 
with stable inverse $-\widetilde{\nu}_X:X \to BTOP$.\\
{\rm (i)} The (rational) {\it Pontrjagin classes} of $X$
are the Pontrjagin classes of $-\widetilde{\nu}_X$
$$p_k(X)~=~p_k(-\widetilde\nu_X) \in H^{4k}(X;\Q)~~(k \geq 0)~.$$
{\rm (ii)} The {\it ${\mathcal L}$-genus} of $X$ is the
${\mathcal L}$-genus of $-\widetilde{\nu}_X$
$${\mathcal L}(X)~=~{\mathcal L}(-\widetilde\nu_X) \in H^{4*}(X;\Q)~.$$
{\rm (iii)} The {\it ${\mathcal L}^H$-genus} of $X$ is
$${\mathcal L}^H(X)~=~8\,i(X)+{\mathcal L}(X) \in H^{4*}(X;\Q)~,$$
with components
$${\mathcal L}_k^H(X)~=~\begin{cases}
1+8\,i(X) \in H^0(X;\Q)&\text{if $k=0$}\vspace*{2mm}\\
{\mathcal L}_k(X)\in H^{4k}(X;\Q)&\text{if $k\geq 1$~.}
\end{cases}$$}
\end{definition}

The ${\mathcal L}$-genus has the same expression
in terms of the Pontrjagin classes of a homology manifold $X$
as for a topological manifold, with components
$$\begin{array}{l}
{\mathcal L}_0(X)~=~p_0(X)~=~1\in H^0(X;\Q)~,\vspace*{2mm}\\{\mathcal L}_1(X)~=~
\displaystyle{\frac{1}{3}}p_1(X)\in H^4(X;\Q)~,\vspace*{2mm}\\
{\mathcal L}_2(X)~=~
\displaystyle{\frac{1}{45}}(7p_2(X)-p_1(X)^2)\in H^8(X;\Q)~,\text{etc.}
\end{array}$$
The Hirzebruch signature theorem also applies to homology manifolds:

\begin{cor} \label{nineteen}
{\rm (i)} If $X$ is an $m$-dimensional homology manifold and 
$f:M=X_{TOP} \to X$ is the canonical normal map then
$${\mathcal L}(M)~=~f^*{\mathcal L}(X) \in H^{4*}(M;\Q)~.$$
{\rm (ii)} The canonical $\LL^{\bullet}$-orientation of an
$m$-dimensional homology manifold $X$ is given rationally by {\rm (}the
Poincar\'e dual of\/{\rm )} the ${\mathcal L}^H$-genus
$$[X]_{\LL} \otimes 1~=~{\mathcal L}^H(X)
\in H_m(X;\LL^{\bullet})\otimes \Q~=~H_{m-4*}(X;\Q)~=~H^{4*}(X;\Q)~.$$
{\rm {\rm (iii)}} The signature of a $4k$-dimensional homology manifold
$X$ is given by
$$\sigma^*(X)~=~\langle{\mathcal L}_k(X),[X]\rangle  \in L^{4k}(\Z)~=~\Z~.$$
\end{cor}
\begin{proof}
(i) Immediate from $f^*(-\widetilde{\nu}_X)=\tau_M$ (stably).\\
(ii) The canonical $\LL^{\bullet}$-orientation of a topological
manifold $M$ is given rationally by the ordinary ${\mathcal L}$-genus
$$[M]_{\LL} \otimes 1~=~{\mathcal L}(M)
\in H_m(M;\LL^{\bullet})\otimes \Q~=~H_{m-4*}(M;\Q)~=~H^{4*}(M;\Q)~,$$
with ${\mathcal L}^H(M)={\mathcal L}(M)$, $i(M)=0$. If $f:M=X_{TOP} \to X$
is the canonical normal map
$$\begin{array}{ll}
[X]_{\LL}\otimes 1\!\!
&=~f_*[M]_{\LL}\otimes 1 + 8\,i(X)\vspace*{2mm} \\
&=~f_*{\mathcal L}(M) + 8\,i(X)~=~{\mathcal L}(X) + 8\,i(X)~
=~{\mathcal L}^H(X) \in H^{4*}(M;\Q)~.
\end{array}$$
(iii) Immediate from (ii), the
identity $\sigma^*(X)=A([X]_{\LL})$ of Theorem \ref{sixteen}, and
the fact that the simply-connected assembly map
$$\begin{array}{l}
A~:~H_{4k}(X;\LL^{\bullet})~=~
H_{4k}(X;L^0(\Z))\oplus H_{4k}(X;\LL^{\bullet}\langle 1 \rangle)\\[1ex]
\hskip100pt
\xymatrix@C+10pt{\ar[r]^-{(0~p_*)}&} 
H_{4k}(\{*\};\LL^{\bullet}\langle 1 \rangle)~=~L^{4k}(\Z)
\end{array}$$
sends $8\,i(X) \in H_{4k}(X;L^0(\Z))$ to 0, with $p:X \to \{*\}$ the unique map.
\end{proof}

\smallskip
\begin{cor} Let $X=M_i$ be the $m$-dimensional homology manifold
with prescribed resolution obstruction $i(X)=i$ given by the
construction of \ref{construct} applied to an $m$-dimensional topological
manifold $M$ {\rm (}$m \geq 6$\/{\rm )}, with normal map $f:X \to M$.\\
{\rm (i)} The canonical $\LL^{\bullet}$-orientation $[X]_{\LL}\in H_m(X;\LL^{\bullet})$
is such that
$$f_*[X]_{\LL}~=~8\,i+ [M]_{\LL} \in H_m(M;\LL^{\bullet})~.$$
{\rm (ii)} The ${\mathcal L}^H$-genus ${\mathcal L}^H(X) \in H^{4*}(X;\Q)=
H_{m-4*}(X;\Q)$ is such that
$$f_*{\mathcal L}^H(X)~=~8\,i + {\mathcal L}(M) \in H^{4*}(M;\Q)$$
\end{cor}
\begin{proof} (i) For any normal map $F:X \to Y$ of
$m$-dimensional homology manifolds
$$F_*[X]_{\LL}~=~(1+T)[F]_{\LL}+ [Y]_{\LL} \in H_m(Y;\LL^{\bullet})~.$$
For $F=f:X \to Y=M$
$$[F]_{\LL}~=~i \in H_m(M;\overline{\LL}_{\bullet})$$
and $(1+T)[F]_{\LL}=8\,i$.\\
(ii) Immediate from (i) and \ref{nineteen}.
\end{proof}
\medskip

We shall now analyze the difference
$$\widetilde{\nu}_{X\times Y} -\widetilde{\nu}_X\times \widetilde{\nu}_Y~:~
X \times Y \to G/TOP$$
for homology manifolds $X,Y$ using the canonical $\LL^{\bullet}$-orientation 
of homology manifolds and the surgery composition and product formulae\,:

\begin{itemize}
\item[{\rm (i)}] (Ranicki \cite{R1} (4.3)) The composite of normal maps
$f:X \to Y$, $g:Y \to Z$ of $n$-dimensional Poincar\'e
spaces is a normal map $gf:X \to Z$ with surgery obstruction
$$\sigma_*(gf )~=~\sigma_*(f) \oplus \sigma_*(g) \in L_n(\Z[\pi_1(Z)])~.$$
\item[{\rm (ii)}] (\cite{R1} (8.1)) 
The product of a normal map $f:M \to X$ of
$m$-dimensional Poincar\'e spaces and a normal map $g:N \to
Y$ of $n$-dimensional Poincar\'e spaces is a normal map $f
\times g:M \times N \to X \times Y$ of $(m+n)$-dimensional geometric
Poincar\'e spaces with surgery obstruction
$$\begin{array}{ll}
\sigma_*(f\times g)&=~\sigma_*(f)\otimes \sigma_*(g) + \sigma^*(X) \otimes
\sigma_*(g) + \sigma_*(f) \otimes \sigma^*(Y)\vspace*{2mm}\\
&=~\sigma^*(M) \otimes \sigma_*(g) + \sigma_*(f) \otimes \sigma^*(Y)
\in L_{m+n}(\Z[\pi_1(X \times Y)])~.
\end{array}$$
\end{itemize}
The formulae are proved on the chain level, using the Eilenberg-Zilber
theorem.

\begin{theorem} \label{seven}
Let $X$ be an $m$-dimensional homology manifold, and let
$Y$ be an $n$-dimensional homology manifold.\\
{\rm (i)} The product $(m+n)$-dimensional homology manifold $X \times
Y$ has canonical $\LL^{\bullet}$-orientation the product of the
canonical $\LL^{\bullet}$-orientations of $X$ and $Y$
$$[X\times Y]_{\LL}~=~[X]_{\LL} \otimes [Y]_{\LL} \in 
H_{m+n}(X \times Y;\LL^{\bullet})~.$$
{\rm (ii)} Let $f:X_{TOP} \to X$, $g:Y_{TOP} \to Y$ be
normal maps from topological manifolds in the normal bordism classes
determined by the canonical $TOP$ reductions
$\widetilde{\nu}_X,\widetilde{\nu}_Y$ of $\nu_X,\nu_Y$, and let
$$h~:~(X\times Y)_{TOP} \to X_{TOP} \times Y_{TOP}$$
be the normal map of topological manifolds classified by
$$\begin{array}{ll}
[h]_{\LL}\!\!\!&=~i(X) \otimes\widetilde{[Y_{TOP}]}_{\LL}+
\widetilde{[X_{TOP}]}_{\LL} \otimes i(Y)\vspace*{2mm} \\
&\in H_{m+n}(X_{TOP} \times Y_{TOP};\LL_{\bullet})~=~
[X_{TOP} \times Y_{TOP},G/TOP]~.
\end{array}$$

The composite normal map
$$(X\times Y)_{TOP} \stackrel{h}{\longrightarrow} X_{TOP} \times Y_{TOP}
\stackrel{f\times g}{\longrightarrow} X \times Y$$
is in the normal bordism class determined by the canonical $TOP$
reduction $\widetilde{\nu}_{X\times Y}$ of $\nu_{X \times Y}$.\\
{\rm (iii)} The canonical $TOP$ reduction
$$\widetilde{\nu}_{X\times Y}~:~X \times Y \to BTOP$$
of the Spivak normal fibration of $X \times Y$
$$\nu_{X\times Y}~=~\nu_X \times \nu_Y~:~X \times Y \to BG$$
differs from the product $TOP$ reduction
$$\widetilde{\nu}_X \times \widetilde{\nu}_Y~:~X \times Y \to BTOP$$
by the element
$$\begin{array}{ll}
t(\widetilde{\nu}_{X\times Y},\widetilde{\nu}_X \times \widetilde{\nu}_Y)\!\!\!
&=~i(X) \otimes\widetilde{[Y]}_{\LL}+
\widetilde{[X]}_{\LL} \otimes i(Y)\vspace*{2mm} \\
&\in H_{m+n}(X \times Y;\LL_{\bullet})~=~[X \times Y,G/TOP]~.
\end{array}$$
\end{theorem}
\begin{proof}
{\rm (i)} This is just the Eilenberg-Zilber theorem on the level of symmetric
Poincar\'e cycles.\\
{\rm (ii)} By construction, $f$ and $g$ are classified by
$$[f]_{\LL}~=~-i(X) \in H_m(X;\overline{\LL}_{\bullet})~~,~~
[g]_{\LL}~=~-i(Y) \in H_n(Y;\overline{\LL}_{\bullet})~.$$

By the surgery product formula (on the level of quadratic Poincar\'e
cycles) the product normal map
$f\times g: X_{TOP} \times Y_{TOP}\to X \times Y$
is classified by
$$\begin{array}{ll}
[f \times g]_{\LL}\!\!\!&=~
[X]_{\LL} \otimes [g]_{\LL} +
[f]_{\LL}\otimes g_*[Y_{TOP}]_{\LL}\vspace*{2mm}\\
&=~-[X]_{\LL}\otimes i(Y) - i(X) \otimes (1+\widetilde{[Y]}_{\LL})
\in H_{m+n}(X \times Y;\overline{\LL}_{\bullet})~.
\end{array}$$
By the surgery composition formula (on the level of quadratic Poincar\'e
cycles)

$$\begin{array}{ll}
[(f\times g)h]_{\LL}\!\!\!
&=~(f\times g)_*[h]_{\LL} +[f\times g]_{\LL}\vspace*{2mm}\\
&=~i(X) \otimes\widetilde{[Y]}_{\LL}+
\widetilde{[X]}_{\LL} \otimes i(Y) \vspace*{2mm}\\
&\hspace*{25mm} 
-[X]_{\LL}\otimes i(Y) - i(X) \otimes g_*[Y_{TOP}]_{\LL}\vspace*{2mm}\\
&=~-i(X)*i(Y)~=~-i(X\times Y) \vspace*{2mm}\\
&\hspace*{10mm}\in
H_{m+n}(X\times Y;L_0(\Z)) \subseteq H_{m+n}(X\times Y;\overline{\LL}_{\bullet})~,
\end{array}$$
so that $(f\times g)h$ is in the canonical normal bordism class.\\
{\rm (iii)} Immediate from {\rm (ii)}, noting that
$$\begin{array}{ll}
[h]_{\LL}\!\!\!
&=~(t(\widetilde{\nu}_{X\times Y},\widetilde{\nu}_X \times \widetilde{\nu}_Y),0)
\vspace*{2mm}\\
&\in H_{m+n}(X_{TOP} \times Y_{TOP};\LL_{\bullet})~=~
H_{m+n}(X \times Y;\LL_{\bullet}) \oplus K_{m+n}
\end{array}$$
with $K_{m+n}$ the $(m+n)$-dimensional homology kernel of $f \times g$.
\end{proof}

\begin{definition}
For any integers $i,j \in \Z$ let
$$i * j~=~i+j+8\,ij \in \Z~,$$
so that
$$(1+8\,i)(1+8j)~=~1+8\,i*j \in \Z~.$$
\end{definition}

\begin{cor} \label{eighteen}

Let $X,Y$ be homology manifolds.\\
{\rm (i)} The ${\mathcal L}^H$-genus of $X \times Y$ is
$${\mathcal L}^H(X \times Y) ~=~{\mathcal L}^H(X)
\otimes {\mathcal L}^H(Y) \in H^{4*}(X \times Y;\Q)~,$$
with
$$\begin{array}{ll}
{\mathcal L}^H_k(X \times Y)~=~
&\sum\limits_{i+j=k}({\mathcal L}_i(X) \otimes {\mathcal L}_j(Y)) +
{\mathcal L}_k(X)\otimes 8\,i(Y)+8\,i(X) \otimes {\mathcal L}_k(Y)\\
& \in H^{4k}(X\times Y;\Q)~(k \geq 1)~.\end{array}$$
{\rm (ii)} {\rm (Quinn \cite{Q2})}
The resolution obstruction of $X \times Y$ is
$$i(X \times Y) ~=~i(X) \ast i(Y) \in \Z~.$$
\end{cor}
\begin{proof} (i) This is just the rationalization of the product
formula of Theorem  \ref{seven} (i).\\
(ii) This is just the 0-dimensional component of the identity of (i).
\end{proof}

As is well-known, the ${\mathcal L}$-genus determines the Pontrjagin
classes, so the formula in Corollary \ref{eighteen} can be used to
determine the Pontrjagin classes of a product homology manifold $X\times Y$
$$p_k(X \times Y) \in H^{4k}(X\times Y;\Q)~~(k \geq 1)$$
in terms of the Pontrjagin classes $p_*(X)$, $p_*(Y)$
and the resolution obstructions $i(X)$, $i(Y)$.

\begin{example} {\normalfont
If $X,Y$ are 4-dimensional homology manifolds, then the
Pontrjagin classes of the product 8-dimensional homology manifold
$X \times Y$ are given by
$$\begin{array}{l}
p_1(X \times Y)~=~p_1(X) \otimes (1+8\,i(Y))+(1+8\,i(X))\otimes p_1(Y) \in H^4(X \times
Y;\Q)~,\vspace*{2mm}\\ p_2(X \times Y)~=~(1+ \displaystyle{\frac{16}{7}}\,
i(X)*i(Y))p_1(X) \otimes p_1(Y) \in H^8(X \times Y;\Q)~.
\end{array}$$}
\end{example}
\qed

\begin{cor} \label{seventeen}
The homology manifold bordism product  for $m\geq 6, n\geq 6$
$$\Omega_m^H(J) \otimes \Omega^H_n(K) \to \Omega_{m+n}^H(J\times K)~;~
X\otimes Y \mapsto X \times Y$$
corresponds under the isomorphisms of \ref{jbord} to the
topological manifold bordism product
$$\begin{array}{l}
\Omega_m^{TOP}(J)[L_0(\Z)] \otimes \Omega^{TOP}_n(K)[L_0(\Z)]
\to \Omega_{m+n}^{TOP}(J \times K)[L_0(\Z)]~;\vspace*{2mm}\\
\hspace*{25mm} M[i]\otimes N[j]\mapsto  (M_i\times N_j)_{TOP}[i*j]~.
\end{array}
$$
Here, $i=I(X)$, $j=I(Y)$, and $M=X_{TOP}\to X$, $N=Y_{TOP}\to Y$ are
the normal maps from topological manifolds determined by the canonical
$TOP$ reductions of $\nu_X$, $\nu_Y$, and $M_i \to M$, $N_j \to N$ are
the normal maps from homology manifolds with resolution obstructions
$I(M_i)=i$, $I(N_j)=j$ given by Proposition \ref{construct}.  The normal
map $h:(X \times Y)_{TOP} \to M \times N $ of Theorem \ref{seven} {\rm (ii)}
classified by
$$[h]_{\LL}~=~i \otimes \widetilde{[N]}_{\LL} +\widetilde{[M]}_{\LL}\otimes j
\in H_{m+n}(M \times N;\LL_{\bullet})$$
is bordant to the composite normal map
$$F~:~(M_i\times N_j)_{TOP} \to M_i \times N_j \to M \times N~.$$
\end{cor}
\noindent{\bf Proof}
The normal maps $M_i \to M$, $N_j \to N$ are classified by
$$[M_i \to M]_{\LL}~=~i \in H_m(M;\overline{\LL}_{\bullet})~~,~~
[N_j \to N]_{\LL}~=~j \in H_n(N;\overline{\LL}_{\bullet})~.$$
Apply the surgery composition and product formulae
to compute the classifying invariant of the
composite $(M_i\times N_j)_{TOP} \to M_i \times N_j \to M \times N$
$$\begin{array}{l}
[(M_i\times N_j)_{TOP} \to M \times N]_{\LL}\vspace*{2mm}\\
\hspace*{10mm}=~
[(M_i\times N_j)_{TOP} \to M_i \times N_j]+[M_i\times N_j\to M \times N]\vspace*{2mm}\\
\hspace*{10mm}=~-i * j +8\,ij +i \otimes [N]_{\LL} + [M]_{\LL} \otimes j \vspace*{2mm}\\
\hspace*{10mm}=~i \otimes \widetilde{[N]}_{\LL} + \widetilde{[M]}_{\LL}\otimes j
\in H_{m+n}(M \times N;\LL_{\bullet})
\subseteq H_{m+n}(M \times N;\overline{\LL}_{\bullet})~,
\end{array}$$
so that $(M_i\times N_j)_{TOP} \to M \times N$, $(X\times Y)_{TOP} \to
M \times N$ are normal bordant.
\qed
\medskip

\begin{remark} {\normalfont Given an $m$-dimensional topological manifold $M$,
an $n$-dimensional topological manifold $N$, and $i,j \in \Z$ let
$$P~=~(M_i\times N_j)_{TOP}$$
be the $(m+n)$-dimensional topological manifold appearing in
Corollary \ref{seventeen}, and let $F:P \to M\times N$ be the degree 1
normal map. The ${\mathcal L}$-genus of $P$ has components
$${\mathcal L}_k(P)\,=\,F^*\big(({\mathcal L}(M) \otimes {\mathcal L}(N))_k
+8\,i \otimes {\mathcal L}_k(N)+ {\mathcal L}_k(M)\otimes 8j\big)\in H^{4k}(P;\Q)~~(k \geq
1)~.$$}
\end{remark}
\qed

\medskip

\begin{example} {\normalfont Take $m=n=4$, $J=K=\{\text{pt.}\}$ in \ref{seventeen}.
The signature defines an isomorphism
$$\omega_4~:~\Omega^{TOP}_4 \to \Z~;~M \mapsto
\displaystyle{\frac{1}{3}}\langle p_1(M),[M]\rangle~=~\sigma^*(M)$$
such that
$$\omega_4(\C\PP^2)~=~1~.$$
Linear combinations of Pontrjagin numbers define an isomorphism
$$\begin{array}{l}
\omega_8~:~\Omega^{TOP}_8\otimes \Q \to \Q\oplus \Q~;\vspace*{2mm}\\
\hphantom{\omega_8~:~}
P \mapsto (\displaystyle{\frac{1}{9}}\langle 5p_2(P)-2p_1(P)^2,[P]\rangle,
\displaystyle{\frac{1}{5}}\langle -2p_2(P)+p_1(P)^2,[P]\rangle)
\end{array}$$
such that
$$\omega_8(\C\PP^2\times \C\PP^2)~=~(1,0)~~,~~\omega_8(\C\PP^4)~=~(0,1)~,$$
with the signature given by
$$\begin{array}{l}
\sigma^*~:~\Omega^{TOP}_8\otimes \Q \stackrel{\omega_8}{\to} \Q\oplus \Q
\stackrel{(1~1)}{\to} L^8(\Z)\otimes \Q~=~\Q~;\vspace*{2mm}\\
P \mapsto
\sigma^*(P)~=~\displaystyle{\frac{1}{9}}\langle 5p_2(P)-2p_1(P)^2,[P]\rangle+
\displaystyle{\frac{1}{5}}\langle -2p_2(P)+p_1(P)^2,[P]\rangle\vspace*{2mm}\\
\hphantom{P \mapsto \sigma^*(P)~}=~
\displaystyle{\frac{1}{45}}\langle 7p_2(P)-p_1(P)^2,[P]\rangle~.
\end{array}
$$
If $M,N$ are 4-dimensional topological manifolds then
$$\omega_8(M \times N)~=~(\sigma^*(M)\sigma^*(N),0)$$
and if $P=(M_i\times N_j)_{TOP}$ for some $i,j \in \Z$ then
$$\begin{array}{l}
\langle p_1(P)^2,[P]\rangle ~=~2(1+8\,i*j)
\langle p_1(M),[M]\rangle\langle p_1(N),[N]\rangle~,\vspace*{2mm}\\
\langle p_2(P),[P]\rangle ~=~(1+\displaystyle{\frac{16}{7}}i*j)
\langle p_1(M),[M]\rangle\langle p_1(N),[N]\rangle~,\vspace*{2mm}\\
\omega_8(P)~=~((1-\displaystyle{\frac{144}{7}}i*j)\sigma^*(M)\sigma^*(N),
(\displaystyle{\frac{144}{7}}i*j)\sigma^*(M)\sigma^*(N))~.
\end{array}
$$}
\end{example}
\qed

The product structure on homology manifold bordism
$$\Omega^H_m\otimes \Omega_n^H \to
\Omega^H_{m+n}~;~X \otimes Y \mapsto X \times Y$$
therefore does not  correspond to the product structure
$$\Omega^{TOP}_m[L_0(\Z)]\otimes \Omega^{TOP}_n[L_0(\Z)] \to
\Omega^{TOP}_{m+n}[L_0(\Z)]~;~M[i] \otimes N[j] \mapsto (M \times N)[i*j]$$
under the abelian group isomorphism of \ref{jbord}
$$\Omega^H_*~ \cong~ \Omega^{TOP}_*[L_0(\Z)]~.$$
However, as stated in Weinberger \cite{Wei}, for $m\geq 6$ there is an
isomorphism of rings
$$\phi~:~\Omega^H_*\otimes \Q ~\cong~(\Omega^{TOP}_*\otimes \Q)[L_0(\Z)]~.$$
We shall now give an explicit description of the isomorphism $\phi$,
which is {\it different} from the isomorphism $\Omega_m^H \cong 
\Omega_m^{TOP}[L_0(\Z)]$ given by \ref{jbord}.

\medskip

Recall that the rational cobordism groups
$$\Omega^{TOP}_*\otimes \Q~=~\Omega^O_*\otimes \Q$$
are detected by Pontrjagin numbers in dimensions $*\equiv 0(\bmod\,4)$, 
and are 0 in dimensions $*\not\equiv 0(\bmod\, 4)$.  
The Pontrjagin numbers of an $m$-dimensional topological manifold $M$ 
with $m \equiv 0(\bmod\,4)$ are
rational linear combinations of the characteristic $\LLL$-numbers
$$\LLL_I(M)~=~\langle\LLL_{i_1}(M)\LLL_{i_2}(M)\cdots \LLL_{i_k}(M),[M]\rangle
\in \Q~,$$
one for each $k$-tuple $I=(i_1,i_2,\dots,i_k)$ of integers $\geq 1$
with $4(i_1+i_2+\dots+i_k)=m$.  Conversely, the characteristic
$\LLL$-numbers of $M$ are rational linear combinations of the
Pontrjagin numbers.  Thus the characteristic $\LLL$-numbers also
determine the rational cobordism class, with isomorphisms
$$\Omega^{TOP}_m\otimes \Q \to \sum\limits_I \Q~;~
M \mapsto \sum\limits_I\LLL_I(M)~.$$

\begin{definition} 
(i) The {\it $\LLL'$-genus} of an $m$-dimensional homology manifold $X$ is
$$\LLL'(X)~=~\frac{1}{1+8\,i(X)}\LLL^H(X) \in H^{4*}(X;\Q)~,$$
with components
$$\LLL'_k(X)~=~
\begin{cases}
1 \in H^0(X;\Q)&\text{if $k=0$}\\
\displaystyle{\frac{1}{1+8\,i(X)}}\LLL_k(X) \in H^{4k}(X;\Q)&\text{if $k\geq 1$~.}
\end{cases}$$
(ii) The {\it $\LLL'$-characteristic numbers} of an $m$-dimensional homology 
manifold $X$ are
$$\begin{array}{c}
\LLL'_I(X)~=~\langle\LLL'_{i_1}(X)\LLL'_{i_2}(X)\dots \LLL'_{i_k}(X),[X]\rangle
\in \Q \vspace*{2mm}\\
(I=(i_1,i_2,\dots,i_k),m=4(i_1+i_2+\dots+i_k))
\end{array}$$
with $\LLL'_I(X)=0$ if $m \not\equiv 0(\bmod\,4)$.\\
(iii) Define a function 
$$\phi~:~\Omega^H_m\otimes \Q \to (\Omega^{TOP}_m\otimes \Q)[L_0(\Z)]~;~
X \mapsto (M\otimes a)[i(X)]$$
using any $m$-dimensional topological manifold $M$ and $a \in \Q$ such that
$$\LLL'_I(X)~=~\LLL_I(M)\otimes a \in \Q$$
for every $I=(i_1,i_2,\dots,i_k)$ with $4(i_1+i_2+\dots+i_k)=m$.
\end{definition}

\medskip
The morphism of \ref{jbord}
$$\psi~:~\Omega^H_m \to \Omega^{TOP}_m[L_0(\Z)]~;~X \mapsto X_{TOP}[i(X)]$$
(which is an isomorphism for $m \geq 6$) is such that the composite
$$\Omega^H_m\otimes \Q \stackrel{\psi}{\to} 
(\Omega^{TOP}_m\otimes \Q)[L_0(\Z)]~\cong~(\sum\limits_I\Q)[L_0(\Z)]$$
sends an $m$-dimensional homology manifold $X$ to 
$$(\sum\limits_I\LLL_I(X))[i(X)]\in (\sum\limits_I\Q)[L_0(\Z)]~.$$
The composite
$$\Omega^H_m\otimes \Q \stackrel{\phi}{\to} 
(\Omega^{TOP}_m\otimes \Q)[L_0(\Z)]~\cong~(\sum\limits_I\Q)[L_0(\Z)]$$
sends an $m$-dimensional homology manifold $X$ to 
$$(\sum\limits_I \LLL_I(M)\otimes a)[i(X)]~=~(\sum\limits_I\LLL'_I(X))[i(X)]
\in (\sum\limits_I\Q)[L_0(\Z)]~,$$
so that $\phi \neq \psi$. The inverse of $\psi$ for $m\geq 6$ is given by
$$\psi^{-1}~:~\Omega^{TOP}_m[L_0(\Z)]\to \Omega^H_m~;~N[i] \mapsto N_i$$
with $N_i$ the homology manifold with resolution obstruction $i$ 
obtained by the construction of \ref{construct}. The composite
$$(\Omega^{TOP}_m\otimes \Q)[L_0(\Z)]\stackrel{\psi^{-1}}{\to} 
\Omega^H_m\otimes \Q  \stackrel{\phi}{\to} 
(\Omega^{TOP}_m\otimes \Q)[L_0(\Z)]~\cong~(\sum\limits_I\Q)[L_0(\Z)]$$
is given by
$$N[i] \mapsto (\sum\limits_I (\frac{1}{1+8\,i})^k{\mathcal L}_I(N_i))[i]~=~
(\sum\limits_I (\frac{1}{1+8\,i})^k{\mathcal L}_I(N))[i]$$
with $I=(i_1,i_2,\dots,i_k)$, $m=4(i_1+i_2+\dots+i_k)$.
\begin{prop}\label{ratprod}
{\rm (i)} The function
$\phi:\Omega^H_m\otimes \Q \to (\Omega^{TOP}_m\otimes \Q)[L_0(\Z)]$
is a ring morphism which is an isomorphism for $m \geq 6$.\\
{\rm (ii)} The rational homology manifold bordism product
$$(\Omega^H_m\otimes \Q) \otimes (\Omega^H_n \otimes \Q) \to 
\Omega^H_{m+n}\otimes \Q~;~(X \otimes Y) \mapsto (X\times Y)$$
corresponds under $\phi$ to the rational topological manifold bordism product
$$\begin{array}{l}
(\Omega^{TOP}_m\otimes \Q)[L_0(\Z)]\otimes (\Omega^{TOP}_n\otimes
\Q)[L_0(\Z)]  \to (\Omega^{TOP}_{m+n}\otimes \Q)[L_0(\Z)]~;\vspace*{2mm}\\
\hspace*{15mm}
(M\otimes a)[i] \otimes (N\otimes b)[j] \mapsto ((M\times N) \otimes ab)[i*j]~.
\end{array}$$
Here $i=i(X)$, $j=i(Y),$ and $M\otimes a=\phi(X)$, $N\otimes b=\phi(Y)$.
\end{prop}
\begin{proof} We know from \ref{jbord} that the morphism of abelian groups
$$\psi~:~\Omega^H_m\otimes \Q \to (\Omega^{TOP}_m\otimes \Q)[L_0(\Z)]~\cong~
(\sum\limits_I\Q)[L_0(\Z)]~;~X\mapsto (\sum\limits_I\LLL_I(X))[i(X)]$$
is an isomorphism for $m \geq 6$. It follows that the morphism of abelian 
groups
$$\phi~:~\Omega^H_m\otimes \Q \to (\Omega^{TOP}_m\otimes \Q)[L_0(\Z)]~\cong~
(\sum\limits_I\Q)[L_0(\Z)]~;~X \mapsto (\sum\limits_I\LLL'_I(X))[i(X)]$$
is also an isomorphism for $m \geq 6$. It remains to show that $\phi$
preserves the multiplicative structures. By \ref{eighteen} (ii)
$$\phi(X\times Y)~=~(P\otimes c)[i*j]$$
for some element $P\otimes c\in \Omega^{TOP}_{m+n}\otimes \Q$. 
We need to show that
$$P\otimes c~=~(M\times N)\otimes ab\in \Omega^{TOP}_{m+n}\otimes \Q~.$$
From the above we see that $M\otimes a$ and $N\otimes b$ correspond to the
$\LLL$-numbers
$$\LLL_I(M\otimes a)~=~\LLL'_I(X)\text{ and }
\LLL_I(N\otimes b)~=~\LLL'_I(Y) \in \Q \text{ respectively.}$$
We compare these $\LLL$-numbers with those of $P\otimes c$, which are
given by
$$\LLL_I(P\otimes c)~=~\LLL'_I(X\times Y) \in \Q~.$$
We observe that by \ref{eighteen} (i)
$$\begin{array}{ll}
\LLL'(X\times Y)\!\!
&=~\displaystyle{\frac{1}{1+8\,i*j}}\,\LLL^H(X\times Y)\vspace*{2mm}\\
&=~\displaystyle{\frac{1}{1+8\,i*j}}\,\LLL^H(X)\otimes \LLL^H(Y)\vspace*{2mm}\\
&=~\LLL'(X)\otimes \LLL'(Y) \in H^{4*}(X\times Y;\Q)~,
\end{array}$$
and that the $\LLL$-numbers satisfy the following product formula
$$\LLL_I((M\times N)\otimes ab)~=~
\sum\limits_{J,K}\LLL_J(M\otimes a)\LLL_K(N\otimes b)\in \Q~,$$
with $I=(i_1,i_2,\dots,i_l)$, and $J=(j_1,j_2,\dots,j_l)$,
$K=(k_1,k_2,\dots,k_l)$ such that
$$j_1+k_1~=~i_1~,~j_2+k_2~=~i_2~,~\dots~,~j_l+k_l=i_l~.$$
Conversely
$$\begin{array}{ll}
\sum\limits_{J,K}
\LLL_J(M\otimes a)\LLL_K(N\otimes b)\!\!
&=~\sum\limits_{J,K}\LLL'_J(X)\LLL'_K(Y)\vspace*{2mm}\\
&=~\sum\limits_{J,K}\LLL'_J(X)\otimes \LLL'_K(Y)\vspace*{2mm}\\
&=~(\LLL'(X)\otimes\LLL'(Y))_I\vspace*{2mm}\\
&=~\LLL'_I(X\times Y)\vspace*{2mm}\\
&=~\LLL_I(P\otimes c) \in \Q~.
\end{array}$$
\end{proof}

\bibliographystyle{amsplain}

\end{document}